\documentclass[pdflatex,sn-basic]{sn-jnl}

\usepackage{graphicx}%
\usepackage{multirow}%
\usepackage{amsmath,amssymb,amsfonts}%
\usepackage{amsthm}%
\usepackage{mathrsfs}%
\usepackage[title]{appendix}%
\usepackage{xcolor}%
\usepackage{textcomp}%
\usepackage{manyfoot}%
\usepackage{booktabs}%
\usepackage{algorithm}%
\usepackage{algorithmicx}%
\usepackage{algpseudocode}%
\usepackage{listings}%
\usepackage{hyperref}
\usepackage{indentfirst}
\theoremstyle{thmstyleone}%
\theoremstyle{thmstyletwo}%

\theoremstyle{thmstylethree}%

\usepackage{tabularx}
\usepackage{amsmath}
\numberwithin{equation}{section}

\begin{document}

\title[Article Title]{Infectious Disease Horizon in Emerging Infectious Diseases}


\author[1,2,3,4]{\fnm{Miaomiao} \sur{Lv}}

\author[1,2,3,4]{\fnm{Juan} \sur{Zhang}}

\author*[1,2,3,4]{\fnm{Zhen} \sur{Jin}} \email{jinzhn@263.net}
\affil[1]{\orgdiv{School of Mathematics and Statistics}, \orgname{Shanxi University}, \orgaddress{\city{Taiyuan 030006}, \state{Shanxi}, \country{China}}}

\affil[2]{\orgdiv{Complex Systems Research Center}, \orgname{Shanxi University}, \orgaddress{\city{Taiyuan 030006}, \state{Shanxi}, \country{China}}}

\affil[3]{\orgdiv{Key Laboratory of Complex Systems and Data Science of Ministry of Education},  \orgaddress{\city{Taiyuan 030006}, \state{Shanxi}, \country{China}}}

\affil[4]{\orgdiv{Shanxi Key Laboratory for Mathematical Technology in Complex Systems},  \orgaddress{\city{Taiyuan 030006}, \state{Shanxi}, \country{China}}}


\abstract{Emerging infectious diseases often spread freely before an outbreak is detected, yet the duration of this free-transmission period and the infection size accumulated during this period are difficult to quantify. These quantities directly affect forecasts of subsequent transmission dynamics, estimates of infection size, and decisions regarding prevention and control measures. In this study, we propose and compute the infectious disease horizon (\textbf{IDH})---the time span from when an infectious disease begins spreading to when it is first detected---and examine its antecedents and consequences. First, based on the SEAIR free-transmission model, we derive transmission characteristics, including the generation interval, the serial interval, the basic reproduction number, and the effective reproduction number, and describe the incidence  using a renewal equation.  Second, by considering two detection routes, we further obtain the distribution and mean of the \textbf{IDH} under different detection settings and quantify the infection size at the time of outbreak detection. Third, intervention measures are incorporated to derive a post-detection transmission model with generally distributed stage durations, followed by the derivation of the corresponding model under exponentially distributed stage durations and the determination of the final infection size. Finally, under the exponential-distribution assumption, the model is applied to the COVID-19 outbreak in Xi'an. The unified framework developed in this study provides a theoretical basis for evaluating early surveillance capacity for emerging infectious diseases, optimizing testing strategies, and determining the intensity of prevention and control measures after outbreak detection.}

\keywords{\textbf{IDH}, Generation interval, Basic reproduction number, Renewal equation, Infection age, Compartmental model}


\pacs[MSC Classification]{92D30, 45D05}

\maketitle

\section{Introduction}\label{sec1}
The frequent occurrence of emerging infectious diseases has imposed a substantial burden on global public health systems and economies \citep{sikkema2026viral}. Early outbreak detection and timely intervention can substantially reduce both the direct and indirect impacts of infectious diseases \citep{bowers2025modelling}. However, the first detected case is not necessarily the individual who was infected earliest, indicating that an outbreak may have already undergone a period of free transmission that is difficult to observe directly before detection \citep{bedford2020cryptic}. Non-pharmaceutical interventions, such as isolation, social distancing, and large-scale screening, are important measures for controlling emerging infectious diseases. The implementation of these measures requires an understanding of transmission before outbreak detection, particularly early characteristics such as the time elapsed from the start of disease transmission to its first detection and the infection size at the time of first detection \citep{liu2023quantitatively,lai2020effect}. In this study, we propose a measure to quantify the time span of undetected transmission and define the time span from when an infectious disease begins spreading to when it is first detected as the infectious disease horizon (\textbf{IDH}).

The \textbf{IDH} is an important temporal characteristic of the early transmission process of emerging infectious diseases but is difficult to observe directly \citep{bowers2025modelling}. Quantifying the \textbf{IDH} can help infer when transmission began, estimate the infection size at first detection, and evaluate the early outbreak detection capacity of existing surveillance systems. It can also provide a quantitative basis for determining screening frequency, optimizing diagnostic strategies, and assessing the appropriate timing for initiating non-pharmaceutical interventions. Therefore, the \textbf{IDH} is an important tool for characterizing the early transmission of emerging infectious diseases and supporting the optimization of surveillance and early warning strategies.

Existing studies have focused primarily on the intervention phase after outbreak detection, whereas relatively less attention has been paid to the free-transmission process before detection \citep{wang2026analysis,jijon2024using}. Studies aiming to infer the onset of an outbreak have often relied on genomic data. However, unsampled transmission chains or lineages that had already gone extinct may cause genomic inference to place the origin of an outbreak later than the actual start of transmission \citep{pekar2022molecular}. In addition, research on the early spread of infectious diseases has generally focused on key epidemiological quantities such as the basic reproduction number, generation interval, serial interval, and growth rate \citep{wallinga2006generation,svensson2007note}. These quantities characterize different aspects of early transmission but do not directly describe the time elapsed before first detection or the infection size accumulated by that time. The generation interval and the basic reproduction number jointly shape the early growth rate \citep{britton2009epidemic}, while the \textbf{IDH}, together with the early growth rate, influences the infection size at first detection, which in turn affects the final infection size. Therefore, linking early transmission characteristics with the \textbf{IDH} and the infection size at first detection can provide a more complete description of the early transmission process of emerging infectious diseases.

During the early stage of an outbreak, when the number of infected individuals is typically small, stochasticity in individual infection processes can substantially affect the generation of secondary infections and the resulting accumulation of infections in the population \citep{van1995stochastic, lloyd2005superspreading, britton2010stochastic, allen2017primer, rebuli2018estimating, czuppon2021stochastic}. Meanwhile, whether and when an infected individual is detected also depend on the route of case detection, with large-scale screening and medical diagnosis of symptomatic infected individuals giving rise to different case-detection processes \citep{wang2026infection, kucirka2020variation}. Therefore, under the joint effects of the transmission process and the case-detection process, the time to first outbreak detection is inherently stochastic, and hence the \textbf{IDH} should be treated as a random variable \citep{bowers2025modelling}. This motivates characterizing the probability distribution and mean of the \textbf{IDH}, as well as analyzing the infection size at first detection and its effects on subsequent transmission \citep{allen2010introduction}.

The remainder of this paper is organized as follows. In Section 2, we characterize the individual infection process and derive key transmission characteristics, including the basic reproduction number, the effective reproduction number, the generation interval, the serial interval, and the incidence. In Section 3, by considering two case-detection routes, we obtain the distribution and mean of the \textbf{IDH} under different detection settings. On this basis, we calculate the infection size at outbreak detection. In Section 4, we derive a post-detection intervention transmission model with generally distributed stage durations and further obtain its formulation under exponential distributions, together with the corresponding final infection size. In Section 5, we use empirical data from the COVID-19 outbreak in Xi'an to estimate the model parameters and analyze the distributional characteristics of the \textbf{IDH}, the infection size at outbreak detection, and the changes in the number of infected individuals before and after detection. We further investigate the effects of the nucleic acid screening interval, the probability of seeking medical care, the infection age of the initially infected individual, the generation interval, and the basic reproduction number on the \textbf{IDH}, as well as how the \textbf{IDH} and the infection size at outbreak detection affect the final infection size. Section 6 summarizes the main findings, discusses the limitations of the study, and outlines directions for future research.

\section{Transmission Characteristics}\label{sec2}
\begin{figure}[H]
	\centering
	\includegraphics[width=0.8\textwidth]{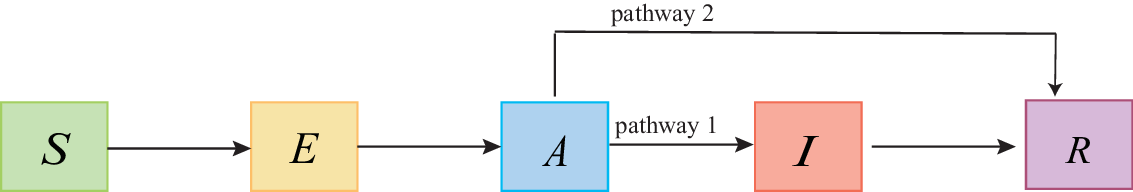}
	\caption{Schematic diagram of individual state transitions. The arrows indicate the possible transition pathways among the susceptible (\(S\)), latent (\(E\)), infectious without symptoms (\(A\)), symptomatic infectious (\(I\)), and recovered (\(R\)) states. After the latent period, an infected individual may follow one of two possible transition pathways, denoted by pathways 1 and 2.}
	\label{model0}
\end{figure}

The model considered in this study is primarily applicable to emerging infectious diseases with a latent period and an infectious period, in which infected individuals may eventually develop symptoms or remain asymptomatic throughout the course of infection. Before an outbreak is first detected, targeted control measures such as isolation and contact tracing have generally not yet been implemented. We therefore refer to this period as the free-transmission stage. According to individual disease states, the population is divided into five states: susceptible (\(S\)), latent (\(E\)), infectious without symptoms (\(A\)), symptomatic infectious (\(I\)), and recovered (\(R\)). After infection, a susceptible individual first enters the latent state \(E\), during which the individual is not infectious. At the end of the latent period, the individual enters the infectious state \(A\) without symptoms. Individuals in state \(A\) subsequently follow one of two pathways. Along the first pathway, the individual develops symptoms with probability \(\rho\), passes through a presymptomatic infectious period, enters the symptomatic infectious state \(I\), and then enters the recovered state \(R\) after losing infectiousness. Along the second pathway, the individual remains asymptomatic throughout the infectious period with probability \(1-\rho\) and enters state \(R\) directly after losing infectiousness. Thus, state \(A\) includes both presymptomatic infectious individuals who will eventually develop symptoms and infectious individuals who remain asymptomatic throughout the infectious period. The transitions among individual disease states are illustrated in Fig.~\ref{model0}. In addition, we consider two routes of case detection: large-scale screening and medical diagnosis of symptomatic infected individuals. Infected individuals without symptoms are mainly detected through large-scale screening conducted at fixed time intervals, whereas symptomatic infected individuals are mainly detected through medical diagnosis.

\begin{figure}[H]
	\centering
	\includegraphics[width=0.8\textwidth]{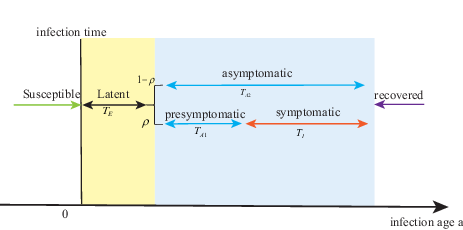}
	\caption{Individual disease course. The yellow region represents the latent period, and the blue region represents the infectious period. After the latent period, an infected individual develops symptoms with probability \(\rho\) after a presymptomatic infectious period, or remains asymptomatic throughout the infectious period with probability \(1-\rho\).}
	\label{individual}
\end{figure}

The infection process of an individual from infection to recovery is illustrated in Fig.~\ref{individual}. In this section, we first define the random variables describing the individual disease course and obtain the probabilities of being in different disease states and the infectivity kernel. We then derive key transmission characteristics, including the generation interval, serial interval, effective reproduction number, basic reproduction number, and incidence, which provide the basis for the subsequent calculation of the \textbf{IDH}.

First, let the random variables $T_E$, $T_{A1}$, $T_{A2}$, and $T_I$ denote the latent period, the time interval from the onset of infectiousness to symptom onset for infected individuals who eventually develop symptoms, the time interval from the onset of infectiousness to recovery for infected individuals who remain asymptomatic throughout the infectious period, and the time interval from symptom onset to recovery for symptomatic infected individuals, respectively. Assume that these stage durations are mutually independent, with probability density functions denoted by $f_E(l)$, $f_{A1}(m)$, $f_{A2}(n)$, and $f_I(r)$, respectively. Their corresponding survival functions are $\overline{F}_E(x)=\Pr(T_E>x), \ 
\overline{F}_{A1}(x)=\Pr(T_{A1}>x), \ 
\overline{F}_{A2}(x)=\Pr(T_{A2}>x), \ 
\overline{F}_I(x)=\Pr(T_I>x)$.
Let $\pi_E(a)$, $\pi_A(a)$, and $\pi_I(a)$ denote the probabilities that an infected individual of infection age $a$ is in the latent state $E$, the infectious state $A$ without symptoms, and the symptomatic infectious state $I$, respectively. The probability of being in state $E$ at infection age $a$ is the probability that the latent period has not yet ended. Thus,
\[
\pi_E(a)
=
\Pr(T_E>a)
=
\overline{F}_E(a).
\]
The probability of being in state $A$ at infection age $a$ consists of two parts: the probability that an individual will eventually develop symptoms but has not yet developed them, and the probability that an individual remains asymptomatic and has not yet recovered. Hence,
\[
\begin{aligned}
	\pi_A(a)
	&=
	\rho\Pr(T_E\leq a<T_E+T_{A1})
	+
	(1-\rho)\Pr(T_E\leq a<T_E+T_{A2})
	\\
	&=
	\rho\int_0^a
	f_E(l)\overline{F}_{A1}(a-l)\,dl
	+
	(1-\rho)\int_0^a
	f_E(l)\overline{F}_{A2}(a-l)\,dl.
\end{aligned}
\]
The probability of being in state $I$ at infection age $a$ is the probability that symptoms have developed but recovery has not yet occurred. Therefore,
\[
\begin{aligned}
	\pi_I(a)
	&=
	\rho\Pr
	\left(
	T_E+T_{A1}\leq a<T_E+T_{A1}+T_I
	\right)
	\\
	&=
	\rho
	\int_0^a
	\int_0^{a-l}
	f_E(l)f_{A1}(m)
	\overline{F}_I(a-l-m)
	\,dm\,dl.
\end{aligned}
\]
Assume that the transmission rate of infectious individuals without symptoms is \(\beta\), whereas that of symptomatic infectious individuals is \(\beta_I=\xi_I\beta\), where \(\xi_I\) is a scaling factor. The  infectivity kernel \(L(a)\) represents the expected infectiousness of an infected individual at infection age \(a\), obtained from the transmission rate at infection age \(a\) weighted by the probability of being infectious at that age \citep{lloyd2009sensitivity}. Based on the probabilities of being in the different infectious states at infection age \(a\), and averaging over the two possible disease-course pathways, the infectivity kernel is given by
\begin{equation}\label{eq:L(a)}
	L(a)=\beta \bigl[\pi_A(a)+\xi_I\pi_I(a)\bigr].
\end{equation}

Next, we use the infectivity kernel to derive the effective reproduction number, the basic reproduction number, and the generation interval. Let \(S^{1}(t)\) denote the number of susceptible individuals at time \(t\). At time \(t\), the rate at which a primary case of infection age \(a\) generates secondary cases is $n(a,t)=S^{1}(t)L(a)$
\citep{nishiura2010time}. Assuming that the current transmission environment remains unchanged, the effective reproduction number, defined as the expected number of secondary infections generated by a single infected individual at time \(t\), is given by
\begin{equation}
	R_t = \int_0^\infty n(a,t)\,da = S^{1}(t)\int_0^\infty L(a)\,da.
	\label{eq:R_t}
\end{equation}
During the early stage of an outbreak, the number of susceptible individuals can be approximated by \(S^{1}(0)\). The basic reproduction number, defined as the average number of secondary infections generated by a single infected individual over the entire infectious period, can be expressed as \citep{diekmann1990definition}
\begin{equation}
	R_0
	=
	S^{1}(0)
	\int_0^\infty L(a)\,da.
	\label{eq:R0}
\end{equation}
The generation interval \(G\) is defined as the time interval between the infection of a primary case and the infection of a secondary case caused by that primary case. The probability density function of the generation interval can be obtained by normalizing the infectivity kernel
\citep{svensson2007note}
\begin{equation}\label{eq:w(a)}
	w_G(a)=\frac{L(a)}{\int_{0}^{\infty} L(a) \, da}= \frac{ n(a,t)}{\int_0^\infty  n(a,t) \, da}=\frac{ n(a,t)}{R_t}, \qquad a\geq 0.
\end{equation}
Clearly, $\int_0^\infty w_G(a)\,da=1$.
Therefore, the mean generation interval is
\begin{equation}
	\overline G
	=
	\int_0^\infty a\,w_G(a)\,da= \frac{\int_{0}^{\infty} a L(a) \, da}{\int_{0}^{\infty} L(a) \, da}.
	\label{eq:mean_generation_interval}
\end{equation}

Since \(T_E, T_{A1}, T_{A2}\), and \(T_I\) are mutually independent, substituting Eq.~(\ref{eq:L(a)}) into Eq.~(\ref{eq:mean_generation_interval}) yields the following explicit expression for the mean generation interval \(\overline{G}\):
\begin{equation}\label{eq:G}
	\overline{G} = \mathbb{E}(T_E) + \frac{ \frac{\rho}{2} \mathbb{E}(T_{A1}^2) + \frac{1 - \rho}{2} \mathbb{E}(T_{A2}^2) +\rho \xi_I \left[\frac{1}{2} \mathbb{E}(T_I^2)+\mathbb{E}(T_{A1}) \mathbb{E}(T_I)\right]}{\rho \mathbb{E}(T_{A1}) + (1 - \rho) \mathbb{E}(T_{A2}) + \rho \xi_I \mathbb{E}(T_I)}.
\end{equation}
Here, \(\mathbb{E}[\cdot]\) denotes mathematical expectation. The mean generation interval consists of two components: the mean latent period and the mean time from the onset of infectiousness to the occurrence of a secondary infection. The latter is obtained by weighting the corresponding times according to the transmission contributions of different disease-course pathways and infectious stages. When \(\rho=1\), all infected individuals eventually develop symptoms and there are no individuals who remain asymptomatic throughout the infectious period. In this case, the mean generation interval reduces to
\[
\overline{G}
=
\mathbb{E}(T_E)
+
\frac{
	\frac{1}{2}\mathbb{E}(T_{A1}^2)
	+
	\xi_I
	\left[
	\frac{1}{2}\mathbb{E}(T_I^2)
	+
	\mathbb{E}(T_{A1})\mathbb{E}(T_I)
	\right]
}{
	\mathbb{E}(T_{A1})
	+
	\xi_I\mathbb{E}(T_I)
}.
\]
When \(\rho=0\), all infected individuals remain asymptomatic throughout the infectious period and do not experience the presymptomatic or symptomatic infectious stages. In this case, the mean generation interval reduces to
\[
\overline{G}
=
\mathbb{E}(T_E)
+
\frac{\mathbb{E}(T_{A2}^{2})}
{2\mathbb{E}(T_{A2})}.
\]
This is consistent with the result in \citep{svensson2007note}.

\begin{figure}[H]
	\centering
	\includegraphics[width=0.8\textwidth]{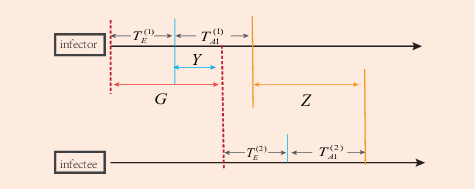}
	\caption{Relationship between the serial interval and the generation interval. The interval between the two red dashed lines represents the generation interval, and the interval between the two orange solid lines represents the serial interval.
	}
	\label{GZ}
\end{figure}

We next discuss the relationship between the generation interval and the serial interval. The serial interval is defined as the time interval between symptom onset in a primary case and symptom onset in the secondary case infected by that primary case \citep{nishiura2010time}. Since the serial interval \(Z\) depends on the symptom onset times of the primary and secondary cases, we consider only transmission pairs in which both cases eventually develop symptoms.
Let $T_E^{(1)}$ and $T_E^{(2)}$ denote the latent periods of the primary and secondary cases, respectively, and let $T_{A1}^{(1)}$ and $T_{A1}^{(2)}$ denote the time intervals from the onset of infectiousness to symptom onset for the primary and secondary cases, respectively. Let $Y$ denote the time from the onset of infectiousness in the primary case to the infection of the secondary case. Then, the generation interval is $G=T_E^{(1)}+Y$.
Taking the infection time of the primary case as the time origin, the symptom onset time of the primary case is $T_E^{(1)}+T_{A1}^{(1)}$,
whereas the symptom onset time of the secondary case is $G+T_E^{(2)}+T_{A1}^{(2)}$.
Therefore, the serial interval $Z$ is
\[
\begin{aligned}
	Z
	&=
	G+
	\left(T_E^{(2)}+T_{A1}^{(2)}\right)
	-
	\left(T_E^{(1)}+T_{A1}^{(1)}\right).
\end{aligned}
\]
The difference between the serial interval and the generation interval is determined by the difference between the incubation periods of the secondary and primary cases. The relationship between the serial interval and the generation interval is illustrated in Fig.\ref{GZ}.

The infectivity kernel, effective reproduction number, basic reproduction number, and generation interval described above are all used to characterize the transmission features associated with the generation of secondary infections. To further describe how the transmission contributions of different infected individuals accumulate in the population and generate new infections, we introduce the incidence and formulate the corresponding renewal equation. Let $\lambda(t)$ denote the incidence at time $t$, that is, the number of new infections per unit time in the population. Suppose that at \(t=0\), a single infected individual who has not yet recovered is introduced into the population, and let \(a_0\) denote the infection age of this individual at the time of introduction. The incidence at time \(t\) consists of contributions from the initially infected individual and all individuals infected during the interval \((0,t]\). Therefore,
\begin{equation}\label{eq:n}
	\lambda(t)=
	S^{1}(t)L(t+a_0)
	+
	S^{1}(t)\int_{0}^{t}\lambda(t-a)L(a)da.
\end{equation}
Substituting Eq.~(\ref{eq:w(a)}) into Eq.~(\ref{eq:n}) gives
\begin{equation}\label{eq:w}
	\lambda(t)=
	R_t w_G(t+a_0)
	+
	\int_0^t R_t \lambda(t-a)w_G(a)da.
\end{equation}
During the early stage of an outbreak, when the number of susceptible individuals remains approximately $S^{1}(0)$, we have $R_t\approx R_0$. The renewal equation can then be written as
\begin{equation}
	\lambda(t)=
	R_0 w_G(t+a_0)
	+
	\int_0^t R_0 \lambda(t-a)w_G(a)da.
	\label{eq:renewal_equation_wG}
\end{equation}
The renewal equation links the infectivity kernel and the generation interval to the incidence, providing the basis for the subsequent calculation of case-detection probabilities, the IDH, and the infection size at outbreak detection.

\section{\textbf{IDH}}\label{sec3}
This section considers two routes through which cases may be detected: medical diagnosis of symptomatic infected individuals and large-scale screening. Detection through medical diagnosis depends on whether an infected individual develops symptoms, the timing of symptom onset, whether medical care is sought, and the delay from symptom onset to diagnosis. Large-scale screening is conducted at fixed time intervals, and its detection probability depends on the timing of screening, the infection age of the individual, and test sensitivity.

We define the time at which the initially infected individual enters the population and transmission begins as $t=0$. Because the individual disease course, secondary transmission, and case-detection process are all stochastic, the time of first case detection is also random. Accordingly, the infectious disease horizon, defined as the time span from when an infectious disease begins spreading to when it is first detected, is a random variable, denoted by $T_H$. In this section, we derive the distribution function and mean of the IDH under different detection routes.

\subsection{\textbf{IDH} under medical diagnosis}\label{subsubsec1}
Symptomatic infected individuals may be detected through medical diagnosis after seeking medical care.  $\rho$ denotes the probability that an infected individual develops symptoms, and $\rho_b$ denotes the probability that a symptomatic infected individual seeks medical care and is diagnosed. Let the random variable \(T_W\) denote the time interval from symptom onset to diagnosis for a symptomatic infected individual, with probability density function \(f_W(w)\) and survival function $\overline{F}_W(x)=\Pr(T_W>x)$. Define the random variable \(T_B\) as the time from infection to detection through medical diagnosis, that is,
\[T_B=
\begin{cases}
	T_E+T_{A1}+T_W,
	& \text{if detected through medical diagnosis},\\
	\infty,
	& \text{otherwise}.
\end{cases}
\]
Since $T_E$, $T_{A1}$, and $T_W$ are mutually independent, the probability that an infected individual of infection age $a$ has been detected through medical diagnosis by time $t$ is
\begin{equation}
	\begin{aligned}
		P_{T_B}(a)
		&=
		\Pr(T_B\leq a)\\
		&=
		\rho\rho_b
		\Pr(T_E+T_{A1}+T_W\leq a)\\
		&=
		\rho\rho_b
		\int_0^a
		\int_0^{a-l}
		\int_0^{a-l-m}
		f_E(l)f_{A1}(m)f_W(w)
		\,dw\,dm\,dl .
	\end{aligned}
	\label{eq:FTB}
\end{equation}
Accordingly, the probability that an infected individual of infection age $a$ has not been detected through medical diagnosis by time $t$ is $
1-P_{T_B}(a)$. 

For the initially infected individual, given that the infection age is $a_0$ at $t=0$ and that the individual has not yet been detected through medical diagnosis, the conditional probability of remaining undetected through medical diagnosis up to time $t$ is
\begin{equation}\label{eq:TBa0}
	G_{B,0}(t)=\Pr(T_B>t+a_0\mid T_B>a_0)
	=
	\frac{
		1-P_{T_B}(t+a_0)
	}{
		1-P_{T_B}(a_0)
	}.
\end{equation}

We next calculate the probability that none of the individuals infected during \((0,t]\) has been detected through medical diagnosis by time \(t\). For any $s\in(0,t]$, \(\lambda(s)\) represents the incidence at time \(s\). Assume that new infection events follow a nonhomogeneous Poisson process with intensity $\lambda(s)$ and that the medical diagnosis processes of different infected individuals are mutually independent. By the thinning property of the Poisson process, the number of individuals infected during $(0,t]$ who have been detected through medical diagnosis by time $t$ follows a Poisson distribution with mean
\[
\int_0^t
\lambda(t-a)P_{T_B}(a)\,da.
\]
Therefore, the probability that none of the individuals infected during $(0,t]$ have been detected through medical diagnosis by time $t$ is
\begin{equation}
	\begin{aligned}
		G_{B,\mathrm{new}}(t)
		&=
		\exp\left\{
		-\int_0^t
		\lambda(t-a)P_{T_B}(a)\,da
		\right\}.
	\end{aligned}
	\label{eq:GBnew}
\end{equation} A detailed derivation is provided in Appendix~\ref{app:medical-diagnosis}.

Considering the initially infected individual and all subsequent infected individuals, the probability that none of the infected individuals in the population has been detected through medical diagnosis by time \(t\) is
\begin{equation}
	G_B(t)
	=
	\frac{
		1-P_{T_B}(t+a_0)
	}{
		1-P  _{T_B}(a_0)
	}
	\exp\left\{
	-\int_0^t
	\lambda(t-a)P_{T_B}(a)\,da
	\right\}.
	\label{eq:GB}
\end{equation}

When only medical diagnosis of symptomatic infected individuals is considered as the route of case detection, the distribution function of the IDH is
\begin{equation}
	\widehat{F}_{T_H}(t)=\Pr(T_H\leq t)
	=
	1-G_B(t).
	\label{eq:FTH_medical}
\end{equation}
Therefore, the mean IDH is
\begin{equation}
	\widehat{\mathbb E}[T_H]
	=
	\int_0^\infty G_B(t)\,dt.
	\label{eq:ETH_medical}
\end{equation}

\subsection{\textbf{IDH} under large-scale screening}\label{subsubsec2}
Large-scale screening refers to the proactive testing of a broad target population to identify infected individuals who have not yet been detected. Depending on the characteristics of the infectious disease, different testing methods may be used, including nucleic acid testing, antigen testing, and serological testing \citep{jeong2025modeling, dimech2021standardization}. Taking COVID-19 as an example, this study mainly considers large-scale nucleic acid screening conducted at a fixed time interval \(\Delta t\). Let \(g(a)\) denote the probability that an infected individual tests positive in a single nucleic acid test at infection age \(a\). Test results for the same infected individual at different screening times are assumed to be mutually independent.

The time at which the initially infected individual enters the population is set as \(t=0\). Large-scale nucleic acid screening is conducted at a fixed interval $\Delta t$, with the first screening occurring at time $t_s$, where $0<t_s<\Delta t$. Thus, the screening times are $t_s,t_s+\Delta t,t_s+2\Delta t,\ldots$. Denote the time of the \(j\)th large-scale nucleic acid screening after the introduction of the initially infected individual by
$D_j=t_s+(j-1)\Delta t$.

The initially infected individual has infection age \(a_0\) at \(t=0\). Since nucleic acid screening is conducted only at the discrete times, the probability that the initially infected individual has not been detected through nucleic acid screening changes only at the screening times and remains constant between consecutive screenings. Therefore, for the initially infected individual, we distinguish the following two cases.

\noindent
(1) When \(0\leq t<t_s\), no nucleic acid screening has yet been conducted. Therefore, the probability that the initially infected individual has not been detected through nucleic acid screening by time \(t\) is 1.

\noindent
(2) When \(D_j\leq t<D_{j+1}\), the initially infected individual has undergone exactly \(j\) nucleic acid screenings, at times \(D_1,\ldots,D_j\), where
$j=\left\lfloor\frac{t-t_s}{\Delta t}\right\rfloor+1$
and $\lfloor\cdot\rfloor$ denotes the floor function. At the \(i\)th screening time \(D_i=t_s+(i-1)\Delta t\), the infection age of the initially infected individual is
\[a_0+D_i=a_0+t_s+(i-1)\Delta t,\quad i=1,\ldots,j\]
Hence, the probability of not being detected at the \(i\)th screening is
\[1-g\left(a_0+t_s+(i-1)\Delta t\right).\]
The probability that the initially infected individual remains undetected through all \(j\) screenings is
\[\prod_{i=1}^{j}
\left[
1-g\left(
a_0+t_s+(i-1)\Delta t
\right)
\right],\]
Accordingly, for any \(t\geq0\), the probability that the initially infected individual has not been detected through nucleic acid screening by time \(t\) is
\begin{equation}\label{eq:TAa0}
	G_{A,0}(t)
	=
	\begin{cases}
		1,
		& 0\leq t<t_s,
		\\[2mm]
		\displaystyle
		\prod_{i=1}^{j}
		\left[
		1-g\left(
		a_0+t_s+(i-1)\Delta t
		\right)
		\right],
		& D_j\leq t<D_{j+1},
		\quad j=1,2,\ldots .
	\end{cases}
\end{equation}

We next calculate the probability that none of the individuals infected during \((0,t]\) has been detected through nucleic acid screening by time \(t\). When \(t=D_j\), let \(a\in[0,D_j)\) denote the infection age of an individual at time \(D_j\), so that the individual was infected at time \(D_j-a\). Since an individual can only be detected by nucleic acid screenings conducted at or after the time of infection, the number of effective nucleic acid screenings experienced by this individual can be written as $k_{D_j}(a)
=
\left\lfloor
\frac{a}{\Delta t}
\right\rfloor+1$. 
The infection age of the individual at the \(i\)th effective screening is
\[
a_{i(a)}:=
a-\bigl(k_{D_j}(a)-i\bigr)\Delta t,
\qquad
i=1,2,\ldots,k_{D_j}(a).
\]
Clearly, \[0\leq a_1<\cdots<a_{k_{D_j}(a)}=a.\]
Therefore, the probability that this individual remains undetected through all effective nucleic acid screenings by time \(D_j\) is
\begin{equation}
	Q_A(a;D_j)
	=
	\prod_{i=1}^{k_{D_j}(a)}
	\left[
	1-g(a_{i(a)})
	\right].
	\label{eq:QA}
\end{equation}
Accordingly, the probability that this individual has been detected through nucleic acid screening by time \(D_j\) is
$
1-Q_A(a;D_j).
$
As in the case of medical diagnosis, by the thinning property of the Poisson process, the number of individuals infected during \((0,D_j]\) who have been detected through nucleic acid screening by time \(D_j\) follows a Poisson distribution with mean
\[\int_0^{D_j}
\lambda(D_j-a)
\left[
1-Q_A(a;D_j)
\right]\,da\]
Hence, the probability that none of the individuals infected during \((0,D_j]\) has been detected through nucleic acid screening by time \(D_j\) is
\begin{equation}
	G_{A,\mathrm{new}}(D_j)
	=
	\exp\left\{
	-\int_0^{D_j}
	\lambda(D_j-a)
	\left[
	1-Q_A(a;D_j)
	\right]\,da
	\right\}.
	\label{eq:GAnew}
\end{equation}
Since nucleic acid screening is conducted only at discrete times, we similarly distinguish the following two cases for individuals infected during the interval (0,t].

\noindent
(1) When $0\leq t<t_s$, no nucleic acid screening
has been conducted, so the probability that none of
these individuals has been detected by time $t$ is $1$.

\noindent
(2) When $D_j\leq t<D_{j+1}$, no additional
large-scale nucleic acid screening is conducted between
the $j$th and $(j+1)$th screening times.
Although new infections may still occur during $(D_j,t]$,
these individuals cannot be detected through nucleic
acid screening during this interval.
Therefore, the probability that none of the individuals
infected during $(0,t]$ has been detected through nucleic
acid screening by time $t$ remains equal to its value
at $D_j$.

Combining the above two cases, we obtain
\[
G_{A,\mathrm{new}}(t)
=
\begin{cases}
	1,
	& 0\leq t<t_s,
	\\[2mm]
	G_{A,\mathrm{new}}(D_j),
	& D_j\leq t<D_{j+1},
	\quad j=1,2,\ldots.
\end{cases}
\]

Considering the initially infected individual and all subsequent
infected individuals, the probability that none of the infected
individuals in the population has been detected through nucleic
acid screening by time $t$ is
\begin{equation}
	G_A(t)
	=
	G_{A,0}(t)
	G_{A,\mathrm{new}}(t).
	\label{eq:GA}
\end{equation}

When only large-scale nucleic acid screening is considered as the route of case detection, the distribution function of the IDH is
\begin{equation}
	\tilde{F}_{T_H}(t)=\Pr(T_H\leq t)
	=
	1-G_A(t).
	\label{eq:FTH_screening}
\end{equation}
Therefore, the mean IDH is
\begin{equation}
	\tilde{\mathbb E}[T_H]
	=
	\int_0^\infty G_A(t)\,dt.
	\label{eq:ETH_screening}
\end{equation}

\subsection{\textbf{IDH} under medical diagnosis and large-scale screening}\label{subsubsec3}
When both medical diagnosis and large-scale nucleic acid screening are considered, we assume that, conditional on the infection process, the nucleic acid screening process and the medical diagnosis process are independent.

For the initially infected individual, according to Eqs.~(\ref{eq:TBa0}) and~(\ref{eq:TAa0}), the probabilities of remaining undetected through large-scale nucleic acid screening and through medical diagnosis by time \(t\) are \(G_{A,0}(t)\) and \(G_{B,0}(t)\), respectively. Therefore, the probability that the initially infected individual has not been detected through either route by time \(t\) is
\begin{equation}
	G_{A,0}(t)G_{B,0}(t).
	\label{eq:GH_initial}
\end{equation}

We next consider individuals infected during \((0,t]\). Since nucleic acid screening is conducted only at the discrete screening times, two cases need to be considered.

\noindent
(1) When \(0\leq t<t_s\), no large-scale nucleic acid screening has yet been conducted. Therefore, individuals infected during \((0,t]\) can only be detected through medical diagnosis. By the thinning property of the Poisson process, the probability that none of these individuals has been detected through either route by time \(t\) is
\begin{equation}\label{eq:GH_new1}
	G_{H,\mathrm{new}}(t)
	=
	\exp\left\{
	-\int_0^t
	\lambda(t-a)P_{T_B}(a)\,da
	\right\},
	\qquad
	0\leq t<t_s.
\end{equation}

\noindent
(2) When \(D_j\leq t<D_{j+1}\), \(j=1,2,\ldots\), \(D_j\) is the most recent large-scale nucleic acid screening time up to time \(t\). 
For an individual with infection age
$a$ at time $t$, we distinguish the following
two cases.

\begin{enumerate}
	\item[(a)] If $0<a<t-D_j$, the individual
	was infected after the $j$th large-scale
	nucleic acid screening and has not undergone
	any nucleic acid screening by time $t$.
	Therefore, the probability that the individual
	has not been detected through either route
	by time $t$ is $1-P_{T_B}(a)$.

	\item[(b)] If $t-D_j\leq a\leq t$, the
	individual had already been infected at the
	time of the $j$th nucleic acid screening,
	and the infection age at time $D_j$ was
	$a-t+D_j$.
	Therefore, the probabilities that the individual
	has not been detected through nucleic acid
	screening and through medical diagnosis by
	time $t$ are $Q_A(a-t+D_j;D_j)$
	and $1-P_{T_B}(a)$, respectively. Under the assumption that the two detection processes are independent, the probability that the individual has not been detected through either route by time \(t\) is
	$$
	Q_H(a,t)
	=
	Q_A(a-t+D_j;D_j)
	\left[1-P_{T_B}(a)\right].
	$$
\end{enumerate}
Accordingly, by the thinning property of the Poisson process, when \(D_j\leq t<D_{j+1}\), the probability that none of the individuals infected during \((0,t]\) has been detected through either route by time \(t\) is
\begin{equation}
	\begin{aligned}
		G_{H,\mathrm{new}}(t)
		=
		\exp\Bigg\{
		&-\int_0^{t-D_j}
		\lambda(t-a)P_{T_B}(a)\,da
		-\int_{t-D_j}^{t}
		\lambda(t-a)
		\left[
		1-Q_H(a,t)
		\right]\,da
		\Bigg\}.
	\end{aligned}
	\label{eq:GH_new2}
\end{equation}

Considering the initially infected individual and all subsequent infected individuals, the probability that none of the infected individuals in the population has been detected through either route by time \(t\) is
\begin{equation}
	\begin{aligned}
		G_H(t)
		&=
		G_{A,0}(t)
		G_{B,0}(t)
		G_{H,\mathrm{new}}(t).
	\end{aligned}
	\label{eq:GH}
\end{equation}
When both medical diagnosis and large-scale nucleic acid screening are considered as routes of case detection, the distribution function of the IDH is
\begin{equation}
	F_{T_H}(t)
	=
	\Pr(T_H\leq t)
	=
	1-G_H(t).
	\label{eq:FTH}
\end{equation}
Therefore, the mean IDH is
\begin{equation}
	\mathbb{E}[T_H]
	=
	\int_0^\infty
	G_H(t)\,dt.
	\label{eq:ETH}
\end{equation}

\subsection{Infection size at outbreak detection}\label{subsec4}
Individuals infected before outbreak detection may be in different disease states when the outbreak is detected. The infection size at outbreak detection depends on the incidence, the individual disease course, and the IDH. Let \(t_l\) denote the mean time of first outbreak detection. Since the time at which the initially infected individual enters the population is set as \(t=0\), we have $t_l=\mathbb{E}[T_H]$. Using the incidence \(\lambda(t)\) and the probabilities that individuals of different infection ages are in different disease states, we calculate the numbers of infected individuals in each state at time \(t_l\).

Let \(E^{1}(t_l)\), \(A^{1}(t_l)\), and \(I^{1}(t_l)\) denote the numbers of individuals at time \(t_l\) in the latent state, the infectious state without symptoms, and the symptomatic infectious state, respectively. Since \(\pi_E(a)\), \(\pi_A(a)\), and \(\pi_I(a)\) denote the probabilities that an individual of infection age \(a\) is in states \(E\), \(A\), and \(I\), respectively, and the initially infected individual has infection age \(t_l+a_0\) at time \(t_l\), the numbers of individuals in the different states at time \(t_l\) are given by
\begin{equation}\label{eq:Etl}
	E^{1}(t_l)
	=\pi_E(t_l+a_0)
	+
	\int_0^{t_l}
	\lambda(t_l-a)\pi_E(a) da,
\end{equation}
\begin{equation}\label{eq:Atl}
	A^{1}(t_l)
	=\pi_A(t_l+a_0)
	+
	\int_0^{t_l}
	\lambda(t_l-a)\pi_A(a) da ,
\end{equation}
and
\begin{equation}\label{eq:Itl}
	I^{1}(t_l)
	=\pi_I(t_l+a_0)
	+
	\int_0^{t_l}
	\lambda(t_l-a)\pi_I(a) da.
\end{equation}
Hence, the total number of infected individuals at the mean outbreak detection time is
\[\begin{aligned}
	N_{\mathrm{inf}}(t_l)
	={}&
	\pi_E(t_l+a_0)
	+\pi_A(t_l+a_0)
	+\pi_I(t_l+a_0)
	\\
	&+
	\int_0^{t_l}
	\lambda(t_l-a)
	\left[
	\pi_E(a)
	+\pi_A(a)
	+\pi_I(a)
	\right]
	\,da .
\end{aligned}
\]

\section{Intervention transmission model}\label{sec4}
\begin{figure}[H]
	\centering
	\includegraphics[width=0.8\textwidth]{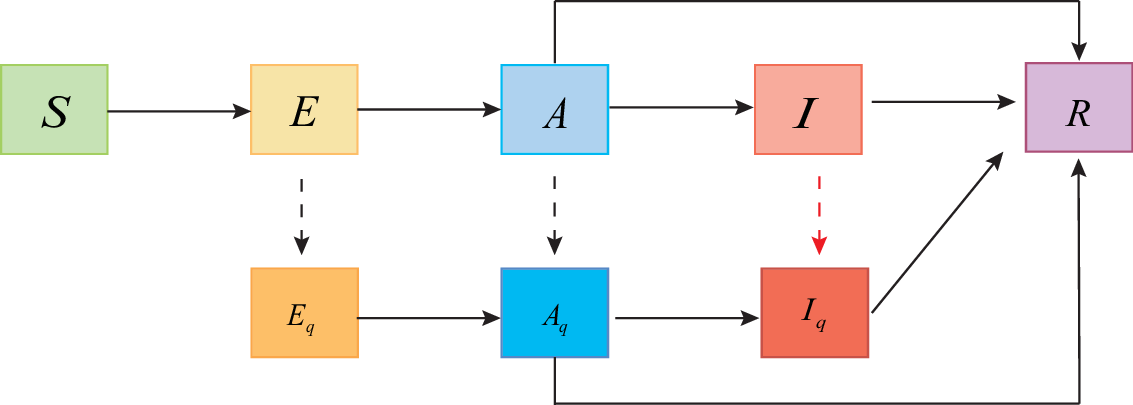}
	\caption{Schematic diagram of individual state transitions after outbreak detection. The black dashed arrows indicate isolation triggered by detection through nucleic acid screening, while the red dashed arrows indicate that symptomatic infected individuals enter isolation after medical diagnosis. The solid arrows represent transitions along the natural disease course.}
	\label{fig:fig4}
\end{figure}

The intervention considered in this section is primarily the isolation of infected individuals. After the outbreak is detected, cases identified through either large-scale nucleic acid screening or medical diagnosis are isolated. Infected individuals in states \(E\) and \(A\) are mainly detected through nucleic acid screening and subsequently enter the corresponding isolated states \(E_q\) and \(A_q\), respectively. Symptomatic infected individuals are mainly detected through medical diagnosis and enter the isolated state \(I_q\) from state \(I\). Isolation prevents further transmission from isolated individuals without altering their natural disease course. Individuals in the isolated latent state \(E_q\) enter the isolated state \(A_q\) after the latent period ends. Those who subsequently develop symptoms enter state \(I_q\) and eventually recover, whereas individuals who remain asymptomatic throughout the infectious period enter state \(R\) directly from \(A_q\). The transitions among individual disease states after outbreak detection are illustrated in Fig.~\ref{fig:fig4}.

\subsection{Under general distributions}\label{subsubsec4}

In this section, we assume that the durations of different disease stages follow general nonnegative continuous distributions and derive the corresponding integral-equation model incorporating isolation measures after outbreak detection.

We reset the time of first outbreak detection to \(t=0\), which serves as the initial time of the post-detection intervention transmission model. At the time of first outbreak detection, all infected individuals are taken to be in the non-isolated states, with isolation initiated immediately thereafter. Accordingly, we set $E_q(0)=A_q(0)=I_q(0)=0$. Let \(E_0\), \(A_0\), and \(I_0\) denote the numbers of infected individuals in states \(E\), \(A\), and \(I\), respectively, at outbreak detection. For analytical convenience, we assume that the individuals already present in states \(E\), \(A\), and \(I\) at outbreak detection have just entered their respective states.

For an individual infected at time $t-a$, let
$Q_{t-a}(b,c)$ denote the probability of remaining
undetected through nucleic acid screening during the
infection-age interval $(b,c]$, where $0\leq b<c\leq a$.
Thus,
\begin{equation}\label{eq:Q}
	Q_{t-a}(b, c) = \prod_{\substack{t-a+b < D_k \leq t-a+c}} [1 - g(D_k - t + a)],
\end{equation}
where  $D_k$ is the time of the $k$th nucleic acid test.
If no nucleic acid screening is conducted during the infection-age interval \((b,c]\), then $Q_{t-a}(b, c)=1$

We first consider individuals infected during \((0,t]\) after outbreak detection and the implementation of intervention measures, and calculate the probabilities that they are in different states at time \(t\). For each state \(X\in\{E,E_q,A,A_q,I,I_q\}\), let \(\Pi_X(t,a)\) denote the probability that an individual with infection age \(a\) at time \(t\) is in state \(X\). Since isolation does not alter the natural disease course, the probability of being in the recovered state is still denoted by \(\Pi_R(a)\).

For a newly infected individual with infection age \(a\) at time \(t\), remaining in the  latent state \(E\) requires two conditions: the latent period has not yet ended, and the individual has not been detected through nucleic acid screening during the infection-age interval \((0,a]\).  \(Q_{t-a}(0,a)\) denotes the probability that an individual remains undetected through nucleic acid screening up to infection age \(a\). Therefore, the probability that the individual is in the latent state \(E\) at time \(t\) with infection age \(a\) is
\begin{equation}
	\Pi_E(t,a)
	=
	\overline F_{E}(a)\,
	Q_{t-a}(0,a).
	\label{eq:piE}
\end{equation}
Correspondingly, the probability that the individual is in the isolated latent state \(E_q\) at time \(t\) with infection age \(a\) is
\begin{equation}
	\Pi_{E_q}(t,a)
	=
	\overline F_{E}(a)
	\left[
	1-Q_{t-a}(0,a)
	\right].
	\label{eq:piEq}
\end{equation}

For an infected individual with infection age \(a\) at time \(t\), remaining in the infectious state \(A\) without symptoms requires that the latent period ends at infection age \(l\in(0,a)\), so that the individual transitions from state \(E\) to state \(A\), and that, by infection age \(a\), the individual has not yet left state \(A\) and has not been detected through nucleic acid screening.  \(Q_{t-a}(0,l)\) and \(Q_{t-a}(l,a)\) denote the probabilities that the individual remains undetected through nucleic acid screening while in states \(E\) and \(A\), respectively. Therefore, the probability that the individual is in the infectious state \(A\) without symptoms at time \(t\) with infection age \(a\) is
\begin{equation}
	\begin{aligned}
		\Pi_A(t,a)
		&=
		\int_0^a
		f_{E}(l)
		\left[
		\rho\,\overline F_{A1}(a-l)
		+
		(1-\rho)\,\overline F_{A2}(a-l)
		\right]\times
		Q_{t-a}(0,l)\,
		Q_{t-a}(l,a)
		\,dl .
	\end{aligned}
	\label{eq:piA}
\end{equation}
Correspondingly, the probability that the individual is in the isolated state \(A_q\) without symptoms at time \(t\) with infection age \(a\) is
\begin{equation}
	\begin{aligned}
		\Pi_{A_q}(t,a)
		&=
		\int_0^a
		f_{E}(l)
		\left[
		\rho\,\overline F_{A1}(a-l)
		+
		(1-\rho)\,\overline F_{A2}(a-l)
		\right]
		\left[
		1-
		Q_{t-a}(0,l)\,
		Q_{t-a}(l,a)
		\right]
		\,dl .
	\end{aligned}
	\label{eq:piAq}
\end{equation}

For an infected individual with infection age \(a\) at time \(t\), remaining in the  symptomatic infectious state \(I\) requires the following conditions. For an individual who eventually develops symptoms, the individual transitions from state \(E\) to state \(A\) at infection age \(l\), remains in state \(A\) for a duration \(m\), and then enters state \(I\). By infection age \(a\), the individual has not yet lost infectiousness, has not been detected through nucleic acid screening before symptom onset, and has not been detected and isolated through medical diagnosis after symptom onset. Therefore, the probability that the individual is in state \(I\) at time \(t\) is
\begin{equation}
	\begin{aligned}
		\Pi_I(t,a)
		&=
		\rho
		\int_0^a
		\int_0^{a-l}
		f_{E}(l)
		f_{A1}(m)
		\overline F_{I}(a-l-m)
		Q_{t-a}(0,l+m)
		\\
		&\qquad\times
		\left[
		(1-\rho_b)
		+
		\rho_b\overline F_{W}(a-l-m)
		\right]
		\,dm\,dl .
	\end{aligned}
	\label{eq:piI}
\end{equation}
Here, $Q_{t-a}(0,l+m)=Q_{t-a}(0,l)\,Q_{t-a}(l,l+m)$, and \(\rho_b\) denotes the probability that a symptomatic infected individual seeks medical care and is diagnosed.
Correspondingly, the probability that the individual is in the isolated symptomatic infectious state \(I_q\) at time \(t\) with infection age \(a\) is
\begin{equation}
	\begin{aligned}
		\Pi_{I_q}(t,a)
		&=
		\rho
		\int_0^a
		\int_0^{a-l}
		f_{E}(l)
		f_{A1}(m)
		\overline F_{I}(a-l-m)
		\\
		&\quad\times
		\left\{
		1-Q_{t-a}(0,l+m)
		+
		\rho_b Q_{t-a}(0,l+m)
		\left[
		1-\overline F_{W}(a-l-m)
		\right]
		\right\}
		\,dm\,dl .
	\end{aligned}
	\label{eq:piIq}
\end{equation}

The probability that an infected individual is in the recovered state \(R\) at infection age \(a\) is
\begin{equation}
	\begin{aligned}
		\Pi_R(a)
		&=
		\rho \Pr(T_E+T_{A1}+T_I\leq a)
		+
		(1-\rho)\Pr(T_E+T_{A2}\leq a)
		\\
		&=
		\rho
		\int_0^a
		\int_0^{a-l}
		\int_0^{a-l-m}
		f_{E}(l)
		f_{A1}(m)
		f_{I}(r)
		\,dr\,dm\,dl
		\\
		&\quad+
		(1-\rho)
		\int_0^a
		\int_0^{a-l}
		f_{E}(l)
		f_{A2}(n)
		\,dn\,dl .
	\end{aligned}
	\label{eq:piR}
\end{equation}
These probabilities satisfy
\[
\Pi_E(t,a)
+
\Pi_{E_q}(t,a)
+
\Pi_A(t,a)
+
\Pi_{A_q}(t,a)
+
\Pi_I(t,a)
+
\Pi_{I_q}(t,a)
+
\Pi_R(a)
=
1.
\]

We next consider infected individuals who are already present at the time of outbreak detection and derive the probabilities that they are in different states at time \(t\). For individuals who are already in state \(E\) at outbreak detection, we approximate them as having just become infected and entered state \(E\) at that time. Therefore, their infection age at time \(t\) after outbreak detection is \(t\). The probabilities that such an individual is in state $
X\in\{E,E_q,A,A_q,I,I_q\}$ at time \(t\) are given by \(\Pi_X(t,t)\), while the probability of being in the recovered state is \(\Pi_R(t)\).

For individuals already in state \(A\) at outbreak detection, we approximate them as having just entered state \(A\) at that time, and denote their infection age upon entering state \(A\) by \(l\). For \(X\in\{A,A_q,I,I_q,R\}\), let \(P_{AX}(t)\) denote the probability that such an individual is in state \(X\) at time \(t\) after outbreak detection. At time \(t\), the infection age of the individual is \(l+t\). Therefore, the probability that the individual remains in the non-isolated state \(A\) at time \(t\) is
\begin{equation}
	P_{AA}(t)
	=
	\left[
	\rho\overline F_{A1}(t)
	+
	(1-\rho)\overline F_{A2}(t)
	\right]
	Q_{t-(l+t)}(l,l+t).
	\label{eq:PAA}
\end{equation}
Similarly, the probability that the individual is in the isolated state \(A_q\) at time \(t\) is
\begin{equation}
	P_{AA_q}(t)
	=
	\left[
	\rho\overline F_{A1}(t)
	+
	(1-\rho)\overline F_{A2}(t)
	\right]
	\left[
	1-Q_{t-(l+t)}(l,l+t)
	\right].
	\label{eq:PAAq}
\end{equation}
The probability that the individual is in the non-isolated symptomatic infectious state \(I\) at time \(t\) is
\begin{equation}
	\begin{aligned}
		P_{AI}(t)
		&=
		\rho
		\int_0^t
		f_{A1}(m)
		Q_{t-(l+t)}(l,l+m)
		\overline F_{I}(t-m)
		\left[
		(1-\rho_b)
		+
		\rho_b\overline F_{W}(t-m)
		\right]
		\,dm .
	\end{aligned}
	\label{eq:PAI}
\end{equation}
The probability that the individual is in the isolated symptomatic infectious state \(I_q\) at time \(t\) is
\begin{equation}
	\begin{aligned}
		P_{AI_q}(t)
		&=
		\rho
		\int_0^t
		f_{A1}(m)
		\overline F_{I}(t-m)
		\Bigg\{
		1-Q_{t-(l+t)}(l,l+m)
		\\
		&\qquad\qquad
		+
		\rho_bQ_{t-(l+t)}(l,l+m)
		\left[
		1-\overline F_{W}(t-m)
		\right]
		\Bigg\}
		\,dm .
	\end{aligned}
	\label{eq:PAIq}
\end{equation}
The probability that the individual is in the recovered state \(R\) at time \(t\) is
\begin{equation}
	P_{AR}(t)
	=\rho
	\int_0^t
	f_{A1}(m)
	\left[
	1-\overline F_{I}(t-m)
	\right]
	\,dm
	+
	(1-\rho)
	\left[
	1-\overline F_{A2}(t)
	\right].
	\label{eq:PAR}
\end{equation}

For individuals who are in state \(I\) at the time of outbreak detection, the probabilities of being in states \(I\), \(I_q\), and \(R\) at time \(t\) are, respectively,
\begin{equation}
	P_{II}(t)
	=
	\overline F_{I}(t)
	\left[
	(1-\rho_b)
	+
	\rho_b\overline F_{W}(t)
	\right],
	\label{eq:PII}
\end{equation}
\begin{equation}
	P_{II_q}(t)
	=
	\rho_b
	\overline F_{I}(t)
	(1-\overline F_{W}(t)),
	\label{eq:PIIq}
\end{equation}
and
\begin{equation}
	P_{IR}(t)
	=
	1-\overline F_{I}(t).
	\label{eq:PIR}
\end{equation}
These probabilities satisfy
\[
P_{AA}(t)+P_{AA_q}(t)
+P_{AI}(t)+P_{AI_q}(t)+P_{AR}(t)
=1,
\]
and
\[P_{II}(t)+P_{II_q}(t)+P_{IR}(t)=1.\]

After outbreak detection, infected individuals in the non-isolated states \(A\) and \(I\) can still infect susceptible individuals, whereas isolated infected individuals no longer contribute to transmission. Therefore, the incidence is
\[
\lambda(t)
=
\beta S(t)
\left[
A(t)+\xi_I I(t)
\right].
\]
Combining the contributions of individuals already present at outbreak detection with those of individuals infected after outbreak detection, we obtain the post-detection intervention transmission model under general distributions
\begin{equation}
	\left\{
	\begin{aligned}
		\frac{dS(t)}{dt}
		&=
		-\lambda(t),
		\\ 
		E(t)
		&=
		E^{(0)}(t)
		+
		\int_0^t
		\lambda(t-a)\Pi_E(t,a)\,da,
		\\ 
		E_q(t)
		&=
		E_q^{(0)}(t)
		+
		\int_0^t
		\lambda(t-a)\Pi_{E_q}(t,a)\,da,
		\\ 
		A(t)
		&=
		A^{(0)}(t)
		+
		\int_0^t
		\lambda(t-a)\Pi_A(t,a)\,da,
		\\ 
		A_q(t)
		&=
		A_q^{(0)}(t)
		+
		\int_0^t
		\lambda(t-a)\Pi_{A_q}(t,a)\,da,
		\\ 
		I(t)
		&=
		I^{(0)}(t)
		+
		\int_0^t
		\lambda(t-a)\Pi_I(t,a)\,da,
		\\ 
		I_q(t)
		&=
		I_q^{(0)}(t)
		+
		\int_0^t
		\lambda(t-a)\Pi_{I_q}(t,a)\,da,
		\\ 
		R(t)
		&=
		R^{(0)}(t)
		+
		\int_0^t
		\lambda(t-a)\Pi_R(a)\,da.
	\end{aligned}
	\right.
	\label{eq:model_general}
\end{equation}
Here, the contributions of individuals already present at outbreak detection to the different states at time \(t\) are
\[
E^{(0)}(t)
=
E_0\Pi_E(t,t),
\
E_q^{(0)}(t)
=
E_0\Pi_{E_q}(t,t),
\]
\[
A^{(0)}(t)
=
E_0\Pi_A(t,t)
+
A_0P_{AA}(t),
\]
\[
A_q^{(0)}(t)
=
E_0\Pi_{A_q}(t,t)
+
A_0P_{AA_q}(t),
\]
\[
I^{(0)}(t)
=
E_0\Pi_I(t,t)
+
A_0P_{AI}(t)
+
I_0P_{II}(t),
\]
\[
I_q^{(0)}(t)
=
E_0\Pi_{I_q}(t,t)
+
A_0P_{AI_q}(t)
+
I_0P_{II_q}(t),
\]
and
\[
R^{(0)}(t)
=
R(0)
+
E_0\Pi_R(t)
+
A_0P_{AR}(t)
+
I_0P_{IR}(t).
\]
Here, \(E_0\), \(A_0\), \(I_0\), and \(R(0)\) denote the numbers of individuals in states \(E\), \(A\), \(I\), and \(R\), respectively, at outbreak detection.

\subsection{Under exponential distributions}\label{subsubsec5}
Assume that the durations of all disease stages and the time
interval from symptom onset to diagnosis follow exponential
distributions, i.e.,
\[
T_E\sim\mathrm{Exp}(\delta),\ 
T_{A1}\sim\mathrm{Exp}(\eta_1),\ 
T_{A2}\sim\mathrm{Exp}(\eta_2),\ 
T_I\sim\mathrm{Exp}(\gamma),\ 
T_W\sim\mathrm{Exp}(\omega).
\]
The corresponding survival functions are
\[
\overline F_E(l)=e^{-\delta l},\ 
\overline F_{A1}(m)=e^{-\eta_1 m},\ 
\overline F_{A2}(n)=e^{-\eta_2 n},\ 
\overline F_I(r)=e^{-\gamma r},\ 
\overline F_{T_W}(w)=e^{-\omega w}.
\]

Since the probability $g(a)$ of testing positive in a nucleic acid test varies with infection age $a$, and the duration of each disease stage is random, individuals within the same disease stage may have different infection ages. For simplification, we take a weighted average of $g(a)$ over infection age to obtain the average probability of testing positive in a single nucleic acid test for individuals in each disease stage. Let $\bar g_E$ and $\bar g_A$ denote the average probabilities of testing positive in a single nucleic acid test while an individual is in the latent state $E$ and the infectious state without symptoms $A$, respectively. Since $\pi_E(a)$ and $\pi_A(a)$ denote the probabilities that an individual of infection age $a$ is in states $E$ and $A$, respectively, under the natural disease course, the average probabilities of testing positive in a single nucleic acid test while in states \(E\) and \(A\) are, respectively,
\[\bar g_E=
\frac{\int_0^\infty g(a)\pi_E(a)\,da}
{\int_0^\infty \pi_E(a)\,da},
\qquad
\bar g_A=
\frac{\int_0^\infty g(a)\pi_A(a)\,da}
{\int_0^\infty \pi_A(a)\,da}.
\]

In addition, after outbreak detection, we apply a
population-level averaging and exponential-equivalent
approximation to the nucleic acid testing process.
Let $\Delta T$ denote the interval between consecutive
nucleic acid tests, with the corresponding testing
frequency $q=1/\Delta T$.  
Because the infection times and testing times of different individuals are not synchronized, testing events are approximated as being uniformly distributed over time at the population level. For an individual who remains in disease state $X$ for a duration $x$, where $X\in\{E,A\}$, the average number of nucleic acid tests received during this interval
is approximated by $x/\Delta T$. Here, $x/\Delta T$ represents the average number of tests at the population level and
therefore need not be an integer.

Under this continuous-time approximation, after replacing the infection age dependent testing positive probability $g(a)$ by the corresponding state-specific average probability
$\bar g_X$, the probability that an individual is not detected through nucleic acid testing while in that state can be approximated by
\[
(1-\bar g_X)^{x/\Delta T}
=
e^{-\kappa_X x},\]
where $
\kappa_X
=
-\frac{\ln(1-\bar g_X)}{\Delta T}.
$
The parameters $\kappa_E$ and $\kappa_A$ represent the equivalent detection rates
obtained by continuous-time approximation of the discrete nucleic acid testing process
at the population level for individuals in the latent state $E$ and the infectious
state without symptoms $A$, respectively.

Therefore, for an individual infected at time $t-a$, if the individual remains in
state $E$ during the infection-age interval $(0,l]$, the probability that the individual
has not yet been detected through nucleic acid screening by infection age $l$ can be
approximated by
\begin{equation}\label{Qa}
	Q_{t-a}(0,l)
	\approx
	(1-\bar g_E)^{l/\Delta T}
	=
	e^{-\kappa_E l}.
\end{equation}
Similarly, if the individual remains in state $A$ during the infection-age interval
$(l,a]$, the probability that the individual has not been detected through nucleic acid
screening by infection age $a$ during this stage can be approximated by
\begin{equation}
	Q_{t-a}(l,a)
	\approx
	(1-\bar g_A)^{(a-l)/\Delta T}
	=
	e^{-\kappa_A(a-l)}.
\end{equation}
Under this population-level averaging and exponential-equivalent
approximation, these probabilities no longer explicitly depend on the
infection time $t-a$.

We now derive the intervention transmission model for the exponential-distribution case based on model~\eqref{eq:model_general}. Unlike the general distribution model with screening at fixed discrete times, under the population-level averaging and exponential-equivalent approximation, the probabilities of remaining undetected through nucleic acid screening no longer explicitly depend on the infection time \(t-a\). Consequently, the state probability kernels depend only on infection age \(a\). For simplicity, we denote the probabilities that an individual of infection age \(a\) is in the different states by \(\Pi_E(a)\), \(\Pi_{E_q}(a)\), \(\Pi_A(a)\), \(\Pi_{A_q}(a)\), \(\Pi_I(a)\), \(\Pi_{I_q}(a)\), and \(\Pi_R(a)\).

Substituting the exponential survival function of \(\overline F_E(a)\) and Eq.~\eqref{Qa} into Eq.~\eqref{eq:piE}, the probability that an individual infected after outbreak detection and the implementation of intervention measures remains in the latent state \(E\) at infection age \(a\) is
\[
\Pi_E(a)
=
e^{-(\delta+\kappa_E)a},
\qquad
\Pi_E(0)=1.
\]
Therefore,
\[
\frac{d\Pi_E(a)}{da}
=
-(\delta+\kappa_E)\Pi_E(a).
\]
Individuals in the latent state \(E\) at time \(t\) arise from two sources: individuals who are already in state \(E\) at outbreak detection and remain in state \(E\) at time \(t\), and individuals infected during \((0,t]\) who are still in state \(E\) at time \(t\). Therefore, the number of infected individuals in the latent state \(E\) at time \(t\) is
\begin{equation}\label{e}
	\begin{aligned}
		E(t)&=E_0\Pi_E(t)+ \int_0^t
		\lambda(t-a)\Pi_E(a)\,da \\
		&= E_0 e^{-(\delta+\kappa_E)t}+\int_0^t \lambda(t-a)e^{-(\delta+\kappa_E)a}\,da.
	\end{aligned}
\end{equation}
Differentiating Eq.~\eqref{e} with respect to \(t\) gives
\begin{equation}
	\begin{aligned}
		\frac{dE(t)}{dt}= \lambda(t)-(\delta+\kappa_E)E(t).
	\end{aligned}\label{de}
\end{equation}
The remaining states \(E_q\), \(A\), \(A_q\), \(I\), \(I_q\), and \(R\) can be derived in the same manner. The detailed derivations are provided in Appendix~\ref{secA2}.

Therefore, when the durations of all disease stages follow exponential distributions, the intervention transmission model~(\ref{eq:model_general}) incorporating isolation measures can be written as the following system:

\begin{equation}
	\left\{
	\begin{aligned}
		\frac{dS(t)}{dt}
		={}&
		-\beta S(t)\left[A(t)+\xi_I I(t)\right],
		\\ 
		\frac{dE(t)}{dt}
		={}&
		\beta S(t)\left[A(t)+\xi_I I(t)\right]
		-(\delta+\kappa_E)E(t),
		\\ 
		\frac{dE_q(t)}{dt}
		={}&
		\kappa_EE(t)-\delta E_q(t),
		\\ 
		\frac{dA(t)}{dt}
		={}&
		\delta E(t)-\kappa_AA(t)
		-\rho\eta_1 K_{1}(t)
		-(1-\rho)\eta_2 K_{2}(t),
		\\ 
		\frac{dA_q(t)}{dt}
		={}&
		\delta E_q(t)+\kappa_AA(t)
		-\rho\eta_1 K_{q1}(t)
		-(1-\rho)\eta_2 K_{q2}(t),
		\\ 
		\frac{dI(t)}{dt}
		={}&
		\rho\eta_1 K_{1}(t)
		-\rho_b\omega\mathcal{J}_{\omega}(t)
		-\gamma I(t),
		\\ 
		\frac{dI_q(t)}{dt}
		={}&
		\rho\eta_1 K_{q1}(t)
		+\rho_b\omega\mathcal{J}_{\omega}(t)
		-\gamma I_q(t),
		\\ 
		\frac{dR(t)}{dt}
		={}&
		\gamma\left[I(t)+I_q(t)\right]
		+
		(1-\rho)\eta_2
		\left[
		K_{2}(t)+K_{q2}(t)
		\right],
	\end{aligned}
	\right.
	\label{eq:jw}
\end{equation}
where
\begin{equation}
	\begin{aligned}
		K_{1}(t)
		={}&
		E_0
		\int_{0}^{t}
		\delta e^{-(\delta+\kappa_E)l}
		e^{-(\eta_1+\kappa_A)(t-l)}
		\,dl
		+
		A_0e^{-(\eta_1+\kappa_A)t}
		\\
		&+
		\int_{0}^{t}
		\lambda(t-a)
		\int_{0}^{a}
		\delta e^{-(\delta+\kappa_E)l}
		e^{-(\eta_1+\kappa_A)(a-l)}
		\,dl\,da,
	\end{aligned}
	\label{eq:K1}
\end{equation}
\begin{equation}
	\begin{aligned}
		K_{2}(t)
		={}&
		E_0
		\int_{0}^{t}
		\delta e^{-(\delta+\kappa_E)l}
		e^{-(\eta_2+\kappa_A)(t-l)}
		\,dl
		+
		A_0e^{-(\eta_2+\kappa_A)t}
		\\
		&+
		\int_{0}^{t}
		\lambda(t-a)
		\int_{0}^{a}
		\delta e^{-(\delta+\kappa_E)l}
		e^{-(\eta_2+\kappa_A)(a-l)}
		\,dl\,da,
	\end{aligned}
	\label{eq:K2}
\end{equation}
\begin{equation}
	\begin{aligned}
		K_{q1}(t)
		={}&
		E_0
		\int_{0}^{t}
		\delta e^{-\delta l}
		e^{-\eta_1(t-l)}
		\left[
		1-e^{-\kappa_E l-\kappa_A(t-l)}
		\right]
		\,dl+
		A_0e^{-\eta_1t}
		\left[
		1-e^{-\kappa_At}
		\right]
		\\
		&+
		\int_{0}^{t}
		\lambda(t-a)
		\int_{0}^{a}
		\delta e^{-\delta l}
		e^{-\eta_1(a-l)}
		\left[
		1-e^{-\kappa_E l-\kappa_A(a-l)}
		\right]
		\,dl\,da,
	\end{aligned}
	\label{eq:Kq1}
\end{equation}
\begin{equation}
	\begin{aligned}
		K_{q2}(t)
		={}&
		E_0
		\int_{0}^{t}
		\delta e^{-\delta l}
		e^{-\eta_2(t-l)}
		\left[
		1-e^{-\kappa_E l-\kappa_A(t-l)}
		\right]
		\,dl+
		A_0e^{-\eta_2t}
		\left[
		1-e^{-\kappa_At}
		\right]
		\\
		&+
		\int_{0}^{t}
		\lambda(t-a)
		\int_{0}^{a}
		\delta e^{-\delta l}
		e^{-\eta_2(a-l)}
		\left[
		1-e^{-\kappa_E l-\kappa_A(a-l)}
		\right]
		\,dl\,da,
	\end{aligned}
	\label{eq:Kq2}
\end{equation}
\begin{equation}
	\begin{aligned}
		\mathcal{J}_{\omega}(t)
		={}&
		\rho E_0
		\int_{0}^{t}
		\int_{0}^{t-l}
		\delta e^{-(\delta+\kappa_E)l}
		\eta_1e^{-(\eta_1+\kappa_A)m}
		e^{-(\gamma+\omega)(t-l-m)}
		\,dm\,dl
		\\
		&+
		\rho A_0
		\int_{0}^{t}
		\eta_1e^{-(\eta_1+\kappa_A)m}
		e^{-(\gamma+\omega)(t-m)}
		\,dm
		+
		I_0e^{-(\gamma+\omega)t}
		\\
		&+
		\rho
		\int_{0}^{t}
		\lambda(t-a)
		\int_{0}^{a}
		\int_{0}^{a-l}
		\delta e^{-(\delta+\kappa_E)l}
		\eta_1e^{-(\eta_1+\kappa_A)m}
		e^{-(\gamma+\omega)(a-l-m)}
		\,dm\,dl\,da.
	\end{aligned}
	\label{eq:Jomega}
\end{equation}
Here, \(\rho K_1(t)\) represents the number of infected individuals who will eventually develop symptoms and are still in state \(A\) at time \(t\), whereas \((1-\rho)K_2(t)\) represents the number of infected individuals who remain asymptomatic throughout the infectious period and are still in state \(A\) at time \(t\). Similarly, \(\rho K_{q1}(t)\) and \((1-\rho)K_{q2}(t)\) represent the corresponding two groups of infected individuals in the isolated state \(A_q\). \(\rho_b\omega\mathcal{J}_{\omega}(t)\) represents the rate at which symptomatic infected individuals transition from state \(I\) to state \(I_q\) through medical diagnosis at time \(t\).

We next derive the ordinary differential equations satisfied by the integral auxiliary functions $K_{1}(t)$, $K_{2}(t)$, $K_{q1}(t)$, $K_{q2}(t)$, and $\mathcal{J}_{\omega}(t)$. First, define the probability kernel corresponding to Eq.~\eqref{eq:K1} as
\[\Pi_{A,1}(a)=
\int_{0}^{a}
\delta e^{-(\delta+\kappa_E)l}
e^{-(\eta_1+\kappa_A)(a-l)}\,dl .\]
Then Eq.~\eqref{eq:K1} can be written as
\[
\begin{aligned}
	K_{1}(t)
	={}&
	E_0\Pi_{A,1}(t)
	+
	A_0e^{-(\eta_1+\kappa_A)t}
	+
	\int_{0}^{t}
	\lambda(t-a)\Pi_{A,1}(a)
	\,da .
\end{aligned}
\]
By the Leibniz rule,
\[
\frac{d\Pi_{A,1}(a)}{da}
=
\delta e^{-(\delta+\kappa_E)a}
-
(\eta_1+\kappa_A)\Pi_{A,1}(a),
\qquad
\Pi_{A,1}(0)=0.
\]
Hence,
\begin{equation}
	\frac{dK_{1}(t)}{dt}
	=
	\delta E(t)
	-
	(\eta_1+\kappa_A)K_{1}(t).
	\label{eq:dk1}
\end{equation}
Similarly, Eqs.~\eqref{eq:K1}--\eqref{eq:Jomega} can be derived using the same approach. The detailed derivations are provided in Appendix~\ref{secA2}.
In summary, under the population-level averaging and exponential-equivalent approximations described above, combining the differential equations satisfied by the five auxiliary integral functions with model~\eqref{eq:jw} yields the following closed system of ordinary differential equations:

\begin{equation}
	\left\{
	\begin{aligned}
		\frac{dS(t)}{dt}
		&=
		-\beta S(t)
		\left[
		\rho K_1(t)
		+
		(1-\rho)K_2(t)
		+
		\xi_I I(t)
		\right],
		\\[1mm]
		\frac{dE(t)}{dt}
		&=
		\beta S(t)
		\left[
		\rho K_1(t)
		+
		(1-\rho)K_2(t)
		+
		\xi_I I(t)
		\right]
		-
		(\delta+\kappa_E)E(t),
		\\ 
		\frac{dE_q(t)}{dt}
		&=
		\kappa_EE(t)-\delta E_q(t),
		\\ 
		\frac{dK_1(t)}{dt}
		&=
		\delta E(t)
		-
		(\eta_1+\kappa_A)K_1(t),
		\\ 
		\frac{dK_2(t)}{dt}
		&=
		\delta E(t)
		-
		(\eta_2+\kappa_A)K_2(t),
		\\ 
		\frac{dK_{q1}(t)}{dt}
		&=
		\delta E_q(t)
		+
		\kappa_AK_1(t)
		-
		\eta_1K_{q1}(t),
		\\ 
		\frac{dK_{q2}(t)}{dt}
		&=
		\delta E_q(t)
		+
		\kappa_AK_2(t)
		-
		\eta_2K_{q2}(t),
		\\ 
		\frac{d\mathcal J_{\omega}(t)}{dt}
		&=
		\rho\eta_1K_1(t)
		-
		(\gamma+\omega)\mathcal J_{\omega}(t),
		\\
		\frac{dI(t)}{dt}
		&=
		\rho\eta_1K_1(t)
		-
		\rho_b\omega\mathcal J_{\omega}(t)
		-
		\gamma I(t),
		\\ 
		\frac{dI_q(t)}{dt}
		&=
		\rho\eta_1K_{q1}(t)
		+
		\rho_b\omega\mathcal J_{\omega}(t)
		-
		\gamma I_q(t),
		\\ 
		\frac{dR(t)}{dt}
		&=
		\gamma
		\left[
		I(t)+I_q(t)
		\right]
		+
		(1-\rho)\eta_2
		\left[
		K_2(t)+K_{q2}(t)
		\right].
	\end{aligned}
	\right.
	\label{eqode}
\end{equation}
The corresponding initial conditions are
\[
S(0)=S_0,\quad
E(0)=E_0,\quad
K_1(0)=K_2(0)=A_0,\quad
I(0)=I_0,
\]
\[
E_q(0)=0,\quad
K_{q1}(0)=K_{q2}(0)=0,\quad
\mathcal J_{\omega}(0)=I_0,\quad
I_q(0)=0.
\]
Here,
\[
\kappa_E
=
-q\ln(1-\bar g_E),
\quad
\kappa_A
=
-q\ln(1-\bar g_A),
\quad
q=\frac{1}{\Delta T}.
\]

\subsection{Final infection size under exponential distributions}\label{subsubsec6}
In this section, we derive the final infection size based on the intervention transmission model~\eqref{eqode} developed above. Let $C$ denote the total number of new infections generated after outbreak detection, which is given by
\begin{equation}\label{eqZ}
	C=S_0-S(\infty).
\end{equation}
Assuming that the outbreak eventually dies out, we have
\[E(\infty)=E_q(\infty)=K_1(\infty)=K_2(\infty)=0,\]
\[K_{q1}(\infty)=K_{q2}(\infty)=\mathcal J_\omega(\infty)=I(\infty)=I_q(\infty)=0.
\]
From the susceptible equation
\[
\frac{dS(t)}{dt}
=
-\beta S(t)
\left[
\rho K_1(t)
+(1-\rho)K_2(t)
+\xi_I I(t)
\right],
\]
dividing both sides by \(S(t)\) and integrating over \([0,\infty)\) gives
\[
\ln\frac{S(\infty)}{S_0}
=
-\beta
\int_0^\infty
\left[
\rho K_1(t)
+(1-\rho)K_2(t)
+\xi_I I(t)
\right]dt.
\]
Therefore,
\begin{equation}
	S(\infty)
	=
	S_0
	\exp\left\{
	-\beta
	\int_0^\infty
	\left[
	\rho K_1(t)
	+(1-\rho)K_2(t)
	+\xi_I I(t)
	\right]dt
	\right\}.
	\label{eqS}
\end{equation}
For notational convenience, define
\[
c_E=\delta+\kappa_E,
\quad
c_{A_1}=\eta_1+\kappa_A,
\quad
c_{A_2}=\eta_2+\kappa_A.
\]
Using
\[
K_1(0)=K_2(0)=A_0,
\quad
\mathcal J_{\omega}(0)=I_0,
\]
and integrating the equations for \(E(t)\), \(K_1(t)\), \(K_2(t)\), \(\mathcal J_{\omega}(t)\), and \(I(t)\) in model~\eqref{eqode} over \([0,\infty)\), we obtain
\begin{equation}
	\int_0^\infty E(t)\,dt
	=
	\frac{C+E_0}{c_E},
	\label{eq:E_inf}
\end{equation}
\begin{equation}
	\int_0^\infty K_1(t)\,dt
	=
	\frac{\delta C}{c_E c_{A_1}}
	+
	\frac{\delta E_0}{c_E c_{A_1}}
	+
	\frac{A_0}{c_{A_1}},
	\label{eq:K1_inf}
\end{equation}
\begin{equation}
	\int_0^\infty K_2(t)\,dt
	=
	\frac{\delta C}{c_E c_{A_2}}
	+
	\frac{\delta E_0}{c_E c_{A_2}}
	+
	\frac{A_0}{c_{A_2}},
	\label{eq:K2_inf}
\end{equation}
\begin{equation}
	\begin{aligned}
		\int_0^\infty I(t)\,dt
		={}&
		\frac{\rho\eta_1}{\gamma}
		\left(
		1-\frac{\rho_b\omega}{\gamma+\omega}
		\right)
		\int_0^\infty K_1(t)\,dt
		+
		\frac{I_0}{\gamma}
		-
		\frac{\rho_b\omega I_0}
		{\gamma(\gamma+\omega)}.
	\end{aligned}
	\label{eq:I_inf}
\end{equation}
Substituting Eqs.~\eqref{eq:E_inf}--\eqref{eq:I_inf}
into Eq.~\eqref{eqS} gives
\begin{equation}\label{Sinfty}
	S(\infty)
	=
	S_0
	\exp\left\{
	-\beta(KC+M)
	\right\},
\end{equation}
where
\[
K
=
\frac{\delta}{c_E}
\left[
\frac{B}{c_{A_1}}
+
\frac{1-\rho}{c_{A_2}}
\right],
\]
and
\[
\begin{aligned}
	M
	={}&
	B
	\left(
	\frac{\delta E_0}{c_Ec_{A_1}}
	+
	\frac{A_0}{c_{A_1}}
	\right)
	+
	(1-\rho)
	\left(
	\frac{\delta E_0}{c_Ec_{A_2}}
	+
	\frac{A_0}{c_{A_2}}
	\right)+
	\frac{\xi_I I_0}{\gamma}
	-
	\frac{
		\xi_I\rho_b\omega I_0
	}{
		\gamma(\gamma+\omega)
	},
\end{aligned}
\]
with
\[
B
=
\rho
\left[
1+
\frac{\xi_I\eta_1}{\gamma}
\left(
1-\frac{\rho_b\omega}{\gamma+\omega}
\right)
\right].
\]
Combining Eq.~\eqref{eqZ},
\[
C=S_0-S(\infty),
\]
with Eq.~\eqref{Sinfty}, we obtain the final size equation satisfied by the
total number of new infections generated after the outbreak is first detected:
\begin{equation}
	C
	=
	S_0
	\left[
	1-
	\exp\left\{
	-\beta(KC+M)
	\right\}
	\right].
	\label{eq:final-size}
\end{equation}
This equation generally does not admit an explicit analytical solution,
and its nonnegative solution $C$ can be obtained numerically. When the transmission rate varies over time, this analytical relation is no longer directly applicable, and the final infection size can be obtained by numerically solving the intervention model.

Since $C$ represents only the number of new infections generated after outbreak detection, the final infection size of the entire outbreak should also include the infected individuals
already present before the outbreak was detected.
Therefore, the final infection size of the entire outbreak is
\begin{equation}
	C_{\mathrm{final}}
	=
	C+E_0+A_0+I_0+R(0).
	\label{eq:total-final-size}
\end{equation}

\section{Applications}\label{sec5}

\subsection{Data and parameter values}\label{subsec6}
In this section, we consider the COVID-19 outbreak caused by the Delta variant in Xi'an from December 15, 2021, to January 20, 2022. The data used in the analysis include the daily numbers of newly reported cases and the permanent resident population. The case data were obtained from the daily epidemic reports released by the Shaanxi Provincial Health Commission \citep{Shaanxi_HC_news_2022}. According to the Xi'an Statistical Yearbook, the permanent resident population of Xi'an in 2021 was 12.873 million \citep{XianStatisticalYearbook2022}. The model parameters were determined based on known epidemiological characteristics of COVID-19, published studies, and the actual prevention and control measures implemented during the Xi'an outbreak.

\begin{figure}[h]
	\centering
	\includegraphics[width=0.55\textwidth]{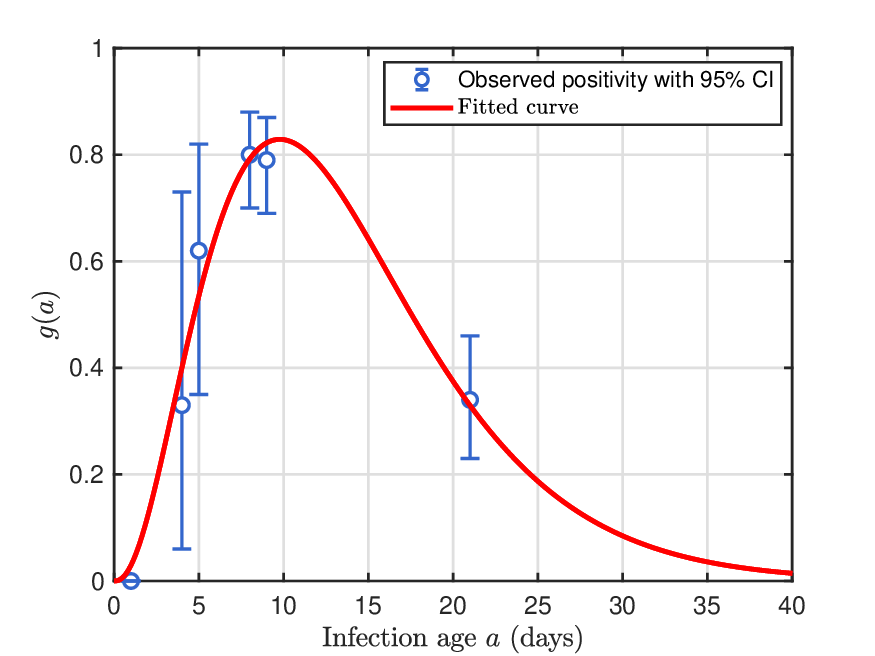}
	\caption{Probability of a positive result in a single nucleic acid test at different infection ages and the corresponding fitted curve.}
	\label{fig:infection age a}
\end{figure}

To characterize how the probability $g(a)$ of testing positive in a single nucleic acid test varies with infection age $a$, we determine $g(a)$ based on previously reported false-negative probabilities of nucleic acid testing at different times after infection. Kucirka et al. estimated the false-negative probability $p_{\mathrm{FN}}(a)$ of RT-PCR testing among infected individuals at different times after exposure \citep{kucirka2020variation}. We approximate the exposure event that resulted in infection as the time of infection and calculate the probability of testing positive in a single nucleic acid test at infection age \(a\) as $g(a)=1-p_{\mathrm{FN}}(a)$. The corresponding values are given in Table~\ref{tab:pcr-positive}. As shown in Table~\ref{tab:pcr-positive}, the probability of a positive nucleic acid test initially increases with infection age and then decreases. To obtain a continuous infection-age-dependent positivity function, we approximate these discrete data using a gamma-type function and estimate \(g(a)\) by least-squares fitting~\citep{johansson2021reducing}, as shown in Fig.~\ref{fig:infection age a}.

\begin{table}[h]
	\caption{False-negative probability and probability of a positive result for a single nucleic acid test at different infection ages~\citep{kucirka2020variation}.}\label{tab:pcr-positive}%
	\begin{tabular}{@{}lll@{}}
		\toprule
		Infection age $a$ (days) &  False-negative probability $p_{\mathrm{FN}}(a)$ &  Positive probability $g(a)$ \\
		\midrule
		1  & 1.00 & 0.00 \\
		4  & 0.67 & 0.33 \\
		5  & 0.38 & 0.62 \\
		8  & 0.20 & 0.80 \\
		9  & 0.21 & 0.79 \\
		21 & 0.66 & 0.34 \\
		\botrule
	\end{tabular}
\end{table}

This wave of the COVID-19 outbreak in Xi'an was mainly caused by the Delta variant. The mean incubation period of individuals infected with the Delta variant was approximately 5.8 days. Meanwhile, the mean latent period from infection to the onset of infectiousness was 3.9 days, and
hence $\delta=\frac{1}{3.9}$~\citep{kang2022transmission}. Therefore, the mean duration from the onset of infectiousness to symptom onset was 1.9 days, giving $\eta_1=\frac{1}{1.9}$. The duration of viral shedding in asymptomatic individuals is generally shorter, and the shedding of infectious virus is concentrated mainly within approximately one week~\citep{murata2021shedding}. Therefore, we
assume $\eta_2=\frac{1}{5}$. The probability that an infected individual eventually develops symptoms is set to $\rho=92\%$~\citep{yu2022proportion}. The mean infectious period of symptomatic individuals is assumed to be 5 days, and hence $\gamma=\frac{1}{5}$. Because transmission after symptom onset is weaker than presymptomatic transmission, we assume $\xi_I=0.1$. Regarding diagnosis through medical care, symptomatic infected individuals may actively seek medical care, with probability \(\rho_{\mathrm{seek}}=35.4\%\)~\citep{liu2021effectiveness}. Accounting for test sensitivity, the probability that an infected individual is diagnosed positive after seeking medical care is assumed to be \(80\%\)~\citep{kucirka2020variation}. Therefore, the probability that a symptomatic infected individual seeks medical care and is diagnosed after symptom onset is \(\rho_b=28.3\%\). The mean time from symptom onset to seeking medical care is 2.3 days, and hence \(\omega=1/2.3\)~\citep{zheng2021health}. Infected individuals who have not yet developed symptoms are generally detected through large-scale nucleic acid screening. After outbreak detection, routine nucleic acid screening was conducted every 3 days, and thus \(\Delta T=3\)~\citep{NationalHealth}.

The actual daily numbers of newly confirmed symptomatic cases are used to estimate the transmission rate parameters in the model. Since control measures were substantially strengthened during the post-detection period, we divide this period into two stages. Xi'an implemented citywide closed management on December 23, 2021. Accordingly, the period from December 15 through December 23 is defined as the first stage, and the second stage begins on December 24. Let \(T_1\) denote the time point separating the two stages. The transmission rate is then assumed to take the following piecewise form:
\[
\beta(t)
=
\begin{cases}
	\beta_1, & 0\leq t<T_1,\\[1mm]
	\beta_2, & t\geq T_1.
\end{cases}
\]
We further assume that the transmission rate during the free-transmission stage before outbreak detection is the same as that during the first stage after detection, and is therefore also given by \(\beta_1\).

Let $\Theta=(\beta_1,\beta_2)$ denote the set of parameters to be estimated. For any candidate
parameter set $\Theta$, we first calculate the mean infectious disease
horizon, $t_l(\Theta)=\mathbb E[T_H(\Theta)]$. Here, $t_l(\Theta)$ is measured from the introduction of the initial infected individual and represents the estimated time of first detection on the time scale of the free-transmission model. The infection states at the time of first detection are then obtained
from the free-transmission model as
\[\begin{aligned}
	S_0= S^{1}\left(t_l;\Theta\right),\
	E_0= E^{1}\left(t_l;\Theta\right),\
	A_0= A^{1}\left(t_l;\Theta\right),\
	I_0= I^{1}\left(t_l;\Theta\right),
\end{aligned}
\]
and $R(0)$ denotes the number of individuals who have already recovered by the time of first detection.

Newly confirmed symptomatic cases arise mainly through two pathways. First, infected individuals who have already been detected and isolated through nucleic acid screening before symptom onset may subsequently develop symptoms and enter state \(I_q\). Second, symptomatic infected individuals who have not yet been isolated may be detected through medical diagnosis and enter state \(I_q\). Let \( C_t\) denote the number of newly reported symptomatic confirmed cases on day \(t\), where \(t=1,\ldots,T\). Accordingly, we fit the model-predicted daily number of individuals newly entering state \(I_q\) to the observed daily number of newly confirmed symptomatic cases. According to model~\eqref{eqode}, the model-predicted number of newly confirmed symptomatic cases on day \(t\) is given by
\[
\widehat C_t(\Theta)
=
\int_{t-1}^{t}
\left[
\rho\eta_1K_{q1}(s;\Theta)
+
\rho_b\omega\mathcal J_{\omega}(s;\Theta)
\right]ds.
\]
The transmission rate parameters are estimated by minimizing the sum of squared residuals between the observed and model-predicted daily numbers of newly confirmed symptomatic cases:
\begin{equation}
	\widehat{\Theta}
	=
	\arg\min_{\Theta}
	\sum_{t=1}^{T}
	\left[
	C_t-\widehat C_t(\Theta)
	\right]^2.
	\label{eq:LS}
\end{equation}
The parameter \(\beta_1\) first determines the pre-detection free-transmission process and thereby affects the infectious disease horizon and the infection size at the time of outbreak detection. The resulting values of \(S_0\), \(E_0\), \(A_0\), \(I_0\), and \(R(0)\) are then used as the initial conditions for the post-detection intervention transmission model. Subsequently, \(\beta_1\) and \(\beta_2\) jointly determine the post-detection transmission dynamics and the daily numbers of newly confirmed symptomatic cases.

\begin{figure}[H]
	\centering
	\includegraphics[width=0.55\textwidth]{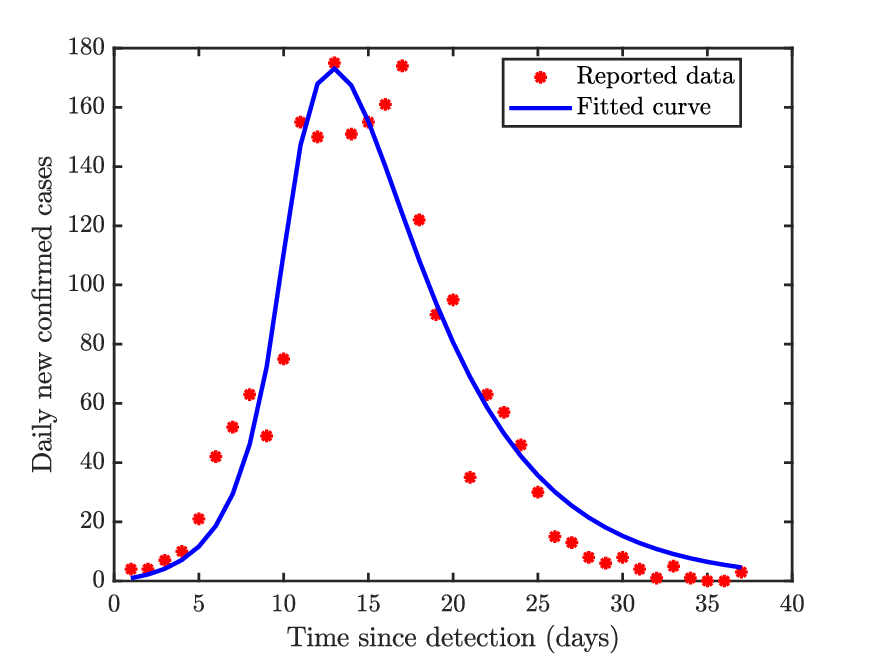}
	\caption{Model fitting of the daily numbers of newly confirmed symptomatic cases during the Xi'an outbreak.
	}
	\label{fig:nihe}
\end{figure}

Therefore, the above procedure provides a unified parameter-estimation framework that links the IDH, the infection size at outbreak detection, and the post-detection case data. The definitions, values, and sources of the model parameters are summarized in Table~\ref{tab:parameter}. The model-fitting results are shown in Fig.~\ref{fig:nihe}. The results indicate that the model can adequately describe the temporal trend in the daily numbers of newly confirmed symptomatic cases during the Xi'an outbreak.

\begin{table}[h]
	\caption{Definitions and values of model parameters.}
	\label{tab:parameter}
	\centering
	\begin{tabular}{@{}lp{7.2cm}ll@{}}
		\toprule
		Parameter & Definition & Value & Source \\
		\midrule
		
		\(\beta_1\) 
		& Transmission rate coefficient in the first stage. 
		& \(2.6917\times10^{-7}\) 
		& Fitted \\
		
		\(\beta_2\) 
		& Transmission rate coefficient in the second stage. 
		& \(1.9118\times10^{-8}\) 
		& Fitted \\
		
		\(\delta\) 
		& Transition rate out of the latent state. 
		& \(1/3.9\) 
		& \citep{kang2022transmission} \\
		
		\(\eta_1\) 
		& Transition rate from the presymptomatic infectious state to symptom onset. 
		& \(1/1.9\) 
		& \citep{kang2022transmission} \\
		
		\(\eta_2\) 
		& Recovery rate of infected individuals who remain asymptomatic throughout the infectious period. 
		& \(1/5\) 
		& Assumed \\
		
		\(\gamma\) 
		& Recovery rate of symptomatic infected individuals. 
		& \(1/5\) 
		& Assumed \\
		
		\(\xi_I\) 
		& Relative infectiousness of symptomatic infected individuals compared with presymptomatic infected individuals. 
		& \(0.1\) 
		& Assumed \\
		
		\(\rho\) 
		& Probability that an infected individual eventually develops symptoms. 
		& \(92\%\) 
		& \citep{yu2022proportion} \\
		
		\(\rho_{\mathrm{seek}}\) 
		& Probability that a symptomatic infected individual seeks medical care.
		& \(35.4\%\) 
		& \citep{liu2021effectiveness} \\
		
		\(\rho_b\) 
		& Probability that a symptomatic infected individual seeks medical care and is diagnosed after symptom onset. 
		& \(28.3\%\) 
		& \citep{liu2021effectiveness,kucirka2020variation} \\
		
		\(\omega\) 
		& Transition rate from symptom onset to diagnosis through medical care. 
		& \(1/2.3\) 
		& \citep{zheng2021health} \\
		\botrule
	\end{tabular}
\end{table}

\subsection{Distributional characteristics of the \textbf{IDH}}\label{subsec7}
\begin{figure}[H]
	\centering
	\includegraphics[width=0.55\textwidth]{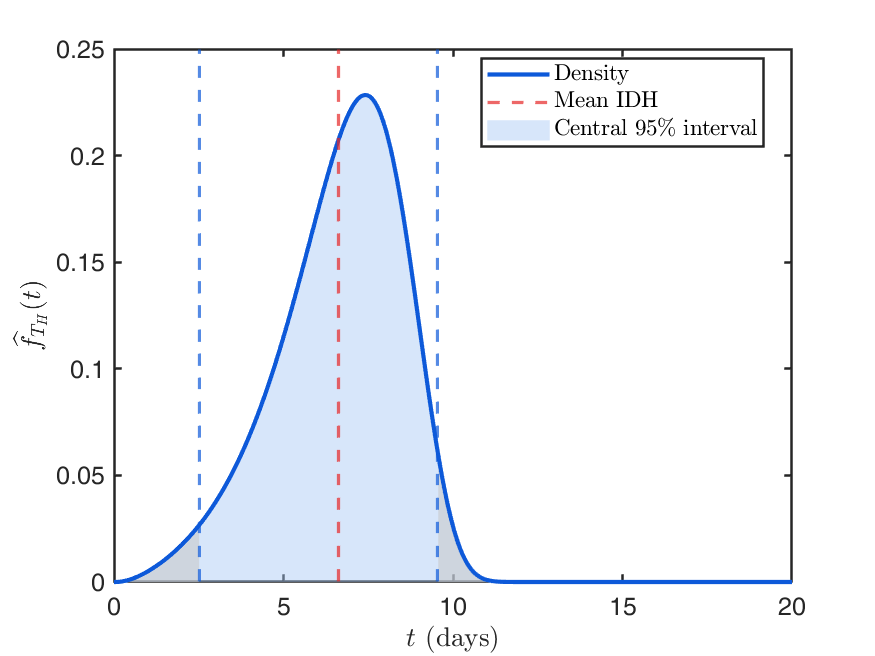}
	\caption{Probability density function \(\hat{f}_{T_H}(t)\) of the IDH. The red dashed line indicates the mean IDH, and the blue shaded region  indicates the central 95\% interval of the IDH distribution.}
	\label{fig:IDH}
\end{figure}

Taking the 2021--2022 COVID-19 outbreak in Xi'an as an example, and assuming that the durations of all disease stages follow exponential distributions, the mean generation interval calculated from Eq.~(\ref{eq:G}) using the parameter values in Table~\ref{tab:parameter} is approximately 7.1 days. 
Since no large-scale nucleic acid
screening targeting this outbreak was conducted before
first detection, we consider medical diagnosis as the
baseline detection route. Fig.~\ref{fig:IDH} shows the IDH distribution,
with a mean of approximately 6.6 days.

\begin{figure}[h]
	\centering
	\includegraphics[width=1\textwidth]{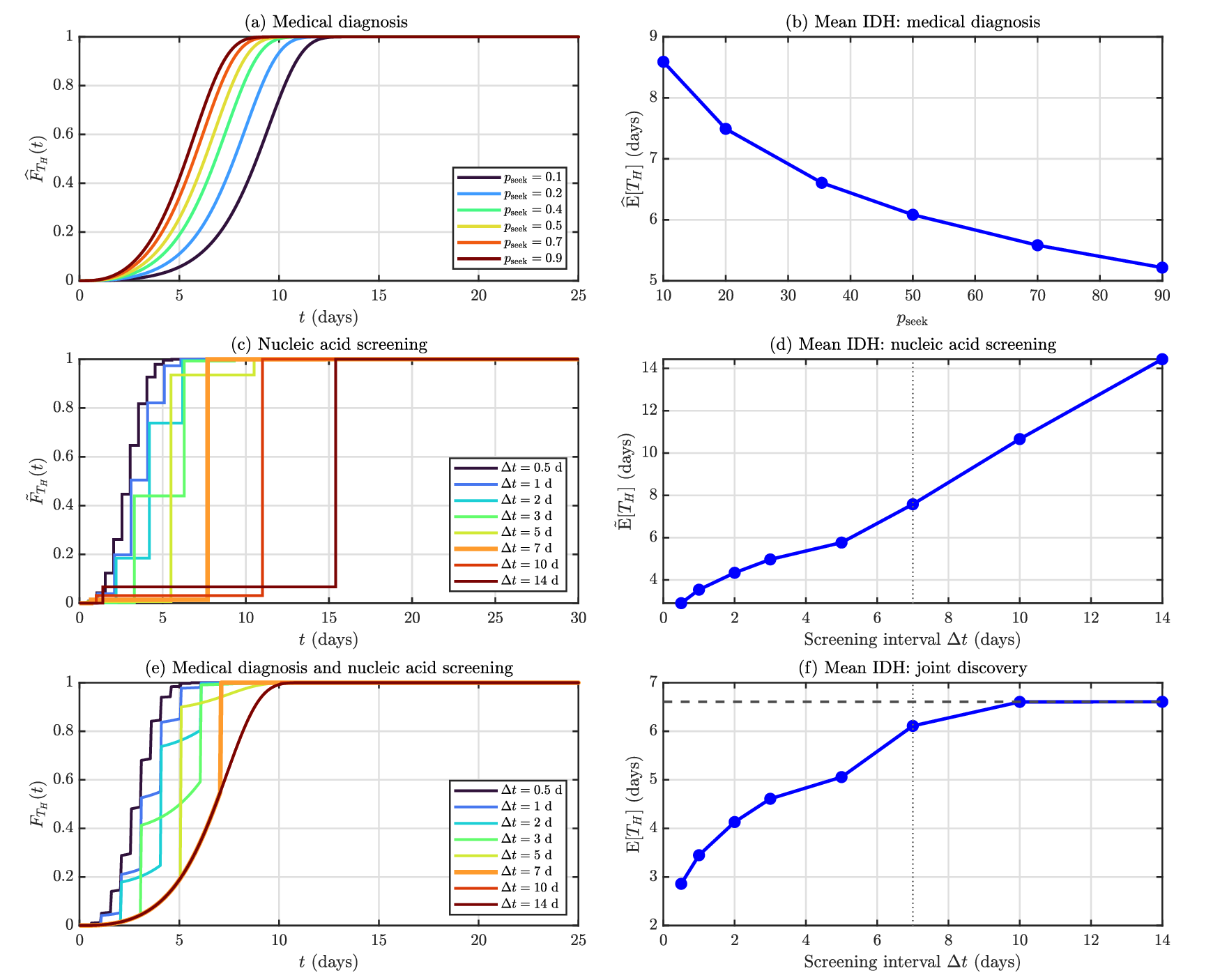}
	\caption{Distribution and mean of the IDH under different detection routes. (a)-(b) Medical diagnosis. (c)-(d) Large-scale nucleic acid screening. (e)-(f) Both medical diagnosis and large-scale nucleic acid screening.
	}
	\label{fig:case-detection processes}
\end{figure}

Large-scale nucleic acid screening and medical diagnosis are the two main routes of outbreak detection considered in this study. Fig.~\ref{fig:case-detection processes} compares the distribution functions and mean  of the IDH under different detection routes. Figs.~\ref{fig:case-detection processes}(a)--(b) show that, under medical diagnosis alone, the mean IDH gradually decreases as the probability that symptomatic infected individuals seek medical care increases. Figs.~\ref{fig:case-detection processes}(c)--(d) show that, under nucleic acid screening alone, a shorter screening interval leads to a smaller mean IDH. When screening is conducted twice per day, the mean IDH is approximately 2.9 days, whereas weekly screening increases it to approximately 7.5 days. Figs.~\ref{fig:case-detection processes}(e)--(f) further show that combining both detection routes shortens the mean IDH compared with either route alone. With the parameters associated with medical diagnosis held fixed, the contribution of nucleic acid screening to outbreak detection gradually weakens as the screening interval increases. When the screening interval reaches two weeks or longer, the mean IDH approaches 6.6 days, corresponding to the case without nucleic acid screening, indicating that outbreak detection depends mainly on medical diagnosis.

\subsection{Infection dynamics}\label{subsec8}
\begin{figure}[H]
	\centering
	
	\begin{minipage}[b]{0.48\textwidth}
		\centering
		\includegraphics[width=\textwidth]{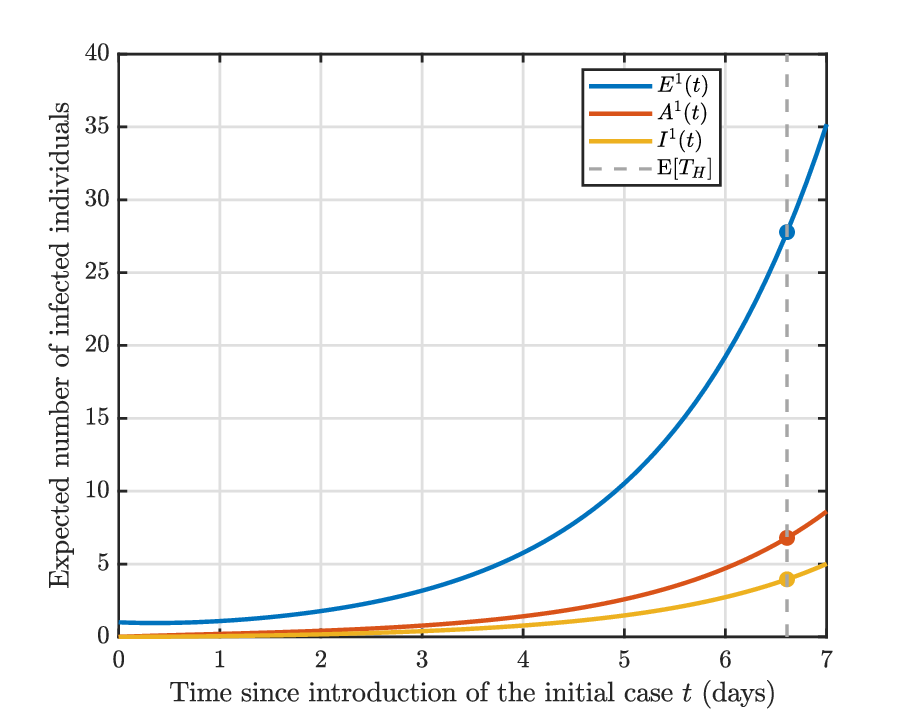}
		
		\small (a)
	\end{minipage}
	\hfill
	\begin{minipage}[b]{0.48\textwidth}
		\centering
		\includegraphics[width=\textwidth]{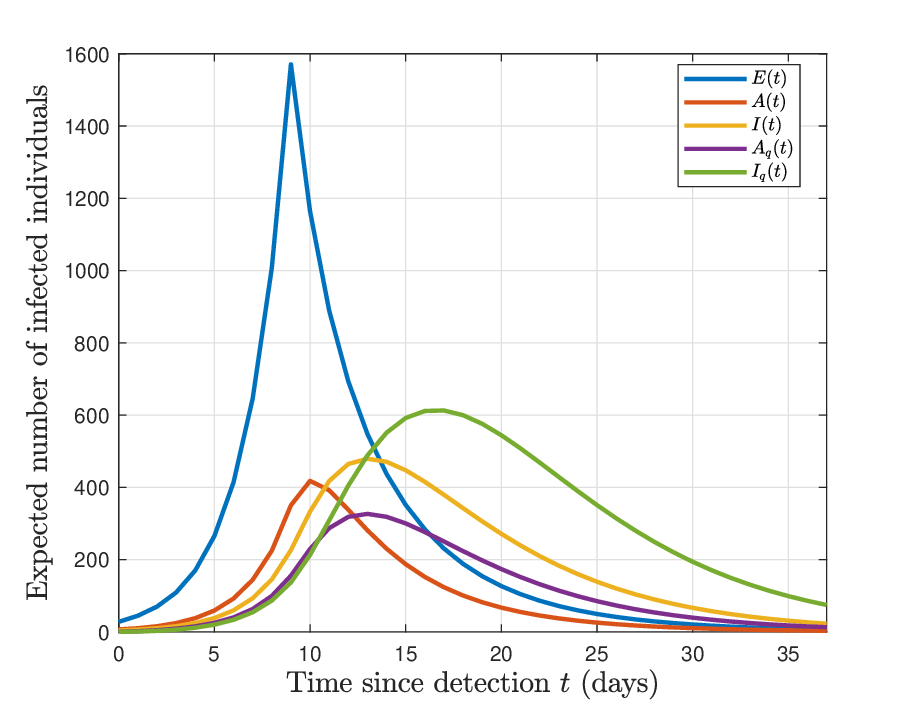}
		
		\small (b)
	\end{minipage}
	
	\caption{Changes in the numbers of infected individuals before and after outbreak detection.
		(a) Changes in the numbers of infected individuals in different infection states during the pre-detection free-transmission stage. The gray dashed line indicates the mean detection time
		$t_l$, and its intersections with the state-specific
		curves indicate the numbers of individuals in states
		$E$, $A$, and $I$ at this time, respectively.
		(b) Changes in the numbers of infected individuals in different infection states after outbreak detection.}
	\label{fig:infection}
\end{figure}

We analyze the infection size at outbreak detection and the transmission dynamics before and after detection. According to Eqs.~(\ref{eq:Etl})--(\ref{eq:Itl}) in Section~\ref{subsec4}, at the mean detection time \(t_l=\mathbb{E}[T_H]\), the numbers of individuals in states \(E\), \(A\), and \(I\) are 28, 7, and 4, respectively, giving a total of 39 infected individuals. Fig.~\ref{fig:infection}(a) shows the temporal changes in the numbers of individuals in different infection states during the pre-detection free-transmission stage. Taking the mean detection time as the connection point between the pre-detection free-transmission model and the post-detection intervention transmission model, we set \(E_0=28\), \(A_0=7\), and \(I_0=4\) as the initial values of the corresponding infection states in the post-detection intervention model. Fig.~\ref{fig:infection}(b) shows the changes in the numbers of individuals in different infection states after outbreak detection.Although intervention measures are implemented after detection, the numbers of infected individuals may continue to increase for a short period before declining as the intervention measures take effect and infected individuals recover.

\begin{figure}[h]
	\centering
	\includegraphics[width=1\textwidth]{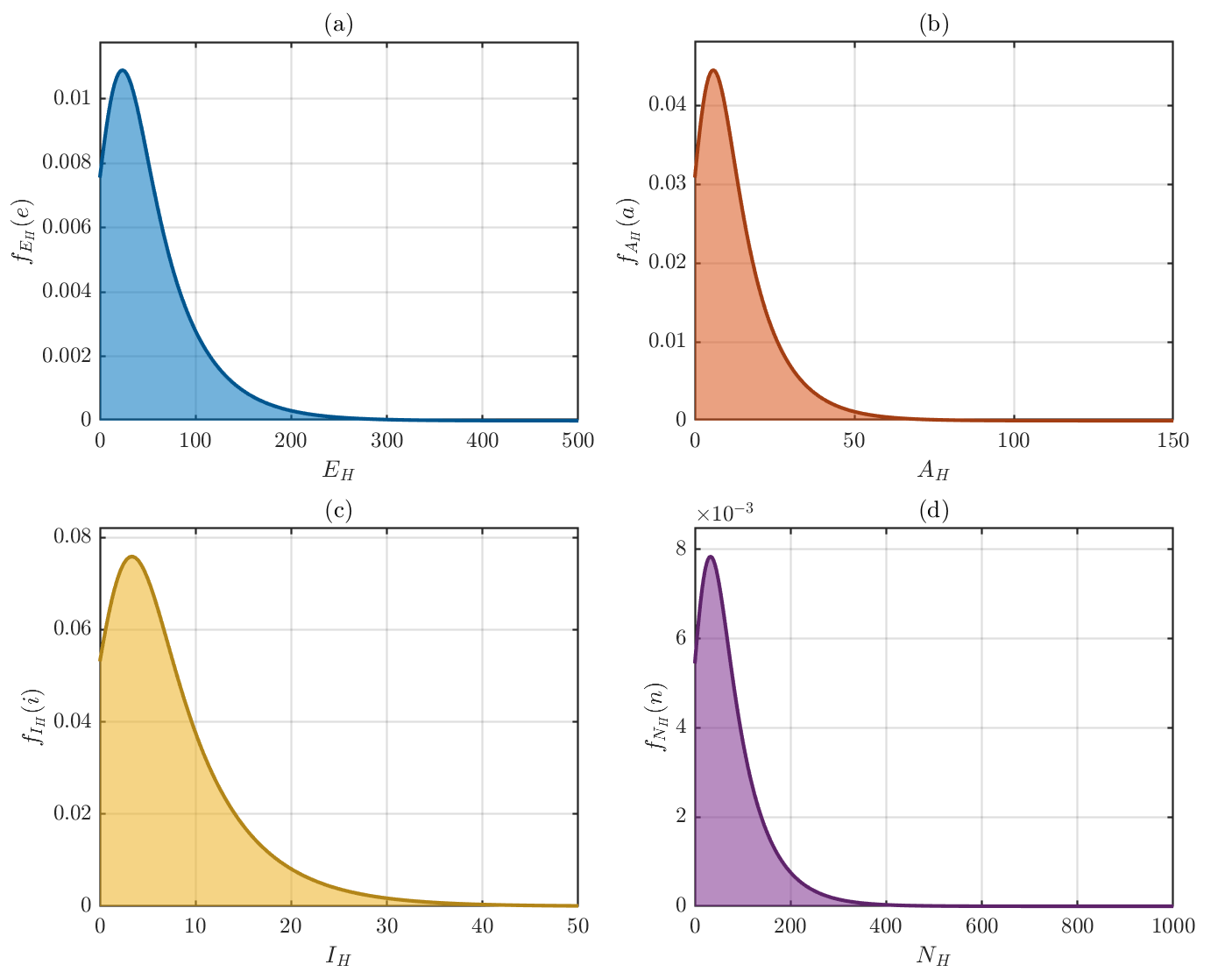}
	\caption{Distributions of the numbers of infected individuals in states \(E\), \(A\), and \(I\), and of the total infection size at outbreak detection.
	}
	\label{fig:state_distributions}
\end{figure}

In addition, the stochasticity of the infectious disease horizon \(T_H\) leads to variation in the time of outbreak detection and consequently to uncertainty in the infection size at detection. Define \(E_H=E^{1}(T_H)\), \(A_H=A^{1}(T_H)\), \(I_H=I^{1}(T_H)\), and \(N_H=E_H+A_H+I_H\). Thus, the distribution of \(T_H\) further induces probability distributions for \(E_H\), \(A_H\), \(I_H\), and \(N_H\) at outbreak detection. Figs.~\ref{fig:state_distributions}(a)-(d) show the distributions of the numbers of infected individuals in the different states and the total infection size at outbreak detection, respectively.

\subsection{Effect of $a_0$ on the mean \textbf{IDH}}\label{subsec9}

\begin{figure}[H]
	\centering
	\includegraphics[width=0.55\textwidth]{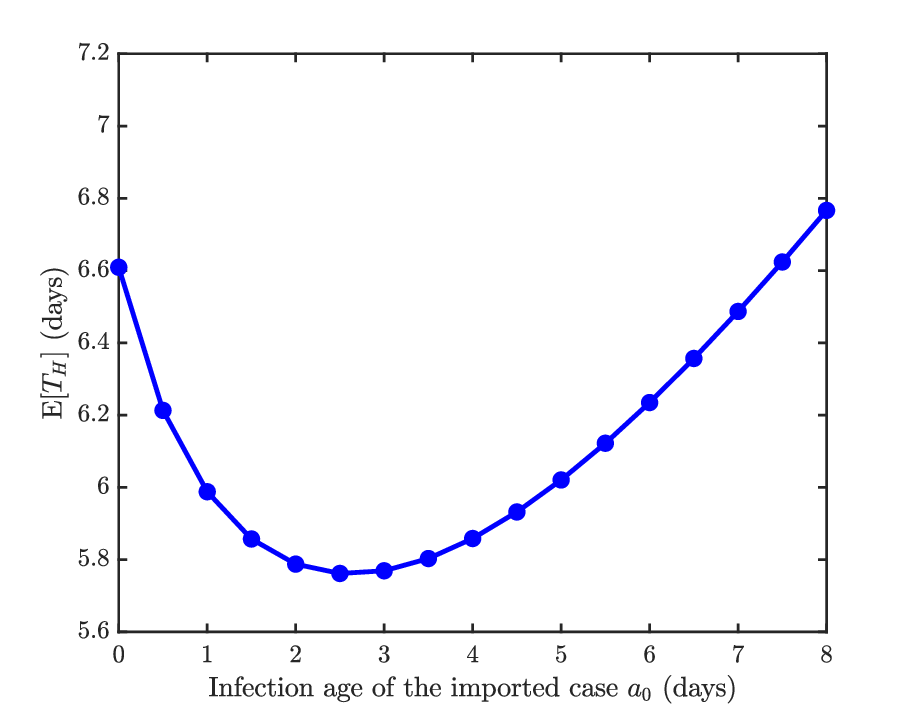}
	\caption{Effect of the infection age of the initial case at introduction on the mean IDH.
	}
	\label{fig:a0}
\end{figure}

The infection age  $a_0$ of the initial case at the time of its introduction affects the mean IDH.
As shown in Fig.~\ref{fig:a0}, \(a_0\) and \(\mathbb{E}[T_H]\) exhibit a non-monotonic relationship. When \(a_0\) is small, the initial case is still in an early stage of infection and requires additional time to reach a stage of higher infectiousness or symptom onset. Consequently, the mean IDH is relatively long. As \(a_0\) increases, the initial case approaches the stage of higher infectiousness and symptom onset, making earlier detection or the generation of secondary infections more likely. The mean IDH therefore decreases and reaches its minimum at approximately \(a_0=3\) days. As \(a_0\) increases further, the initial case gradually passes beyond the main infectious period, and its ability to generate secondary infections declines. Outbreak detection then relies more heavily on subsequent infected individuals, causing the mean IDH to increase again. These results indicate that an infection age that is either too small or too large when the initial case appears may prolong the IDH, whereas the outbreak is more likely to be detected earlier when the initial case is at a stage of relatively high infectiousness and close to symptom onset.

\subsection{Effects of \(\overline G\) and \(R_0\) on the \textbf{IDH}}\label{subsec10}
\begin{figure}[H]
	\centering
	\includegraphics[width=1\textwidth]{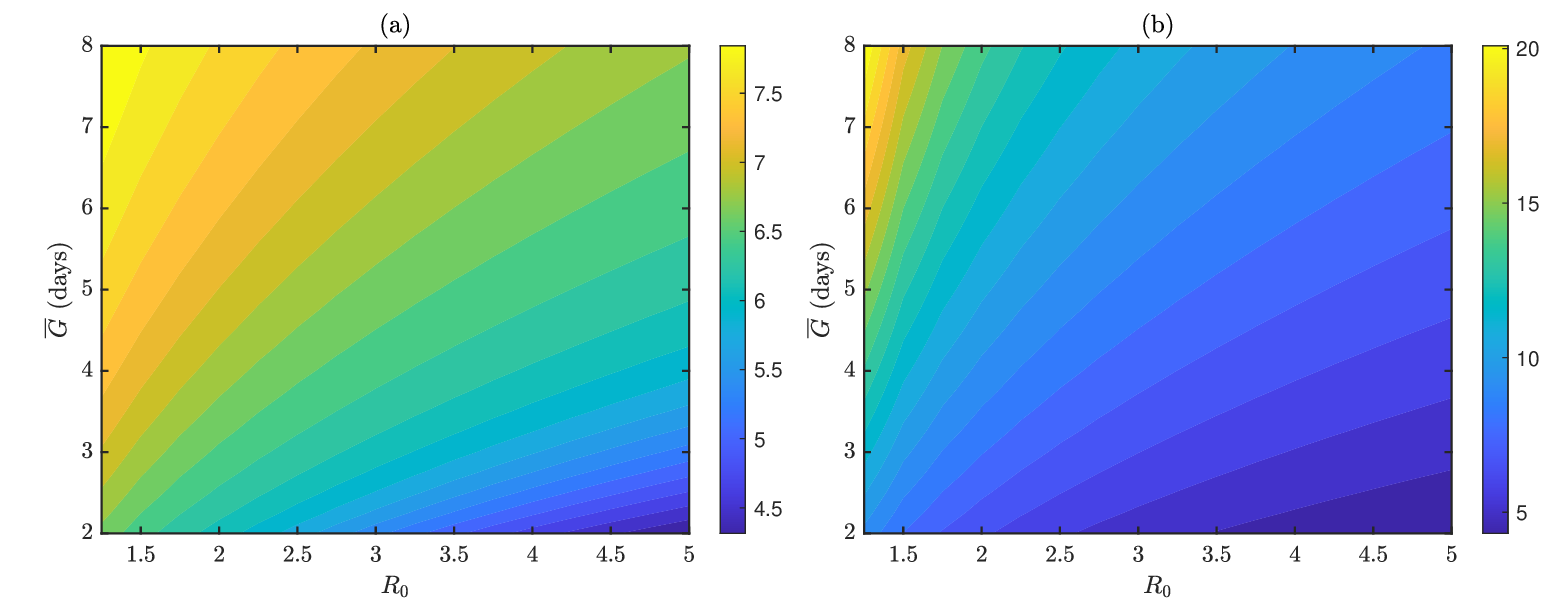}
	\caption{Effects of the mean generation interval and the basic reproduction number on the infectious disease horizon. (a) Medical diagnosis combined with weekly nucleic acid screening. (b) Medical diagnosis only.
	}
	\label{fig:G_R_H}
\end{figure}

Both the basic reproduction number \(R_0\) and the mean generation interval \(\overline G\) affect the IDH. As shown in Fig.~\ref{fig:G_R_H}(a), when weekly nucleic acid screening and medical diagnosis are both implemented, the mean IDH decreases with increasing \(R_0\) for a fixed \(\overline G\), whereas for a fixed \(R_0\), it increases with increasing \(\overline G\). This is because a larger \(R_0\) accelerates the accumulation of infected individuals, while a longer generation interval delays the occurrence of secondary infections. Fig.~\ref{fig:G_R_H}(b) shows that, when outbreak detection relies only on medical diagnosis, the mean IDH is more sensitive to changes in \(R_0\) and \(\overline G\). These results indicate that nucleic acid screening not only shortens the mean IDH but also reduces its sensitivity to changes in transmission characteristics.

\subsection{Effect of the \textbf{IDH} on the final infection size}\label{subsec11}
\begin{figure}[H]
	\centering
	\includegraphics[width=0.55\textwidth]{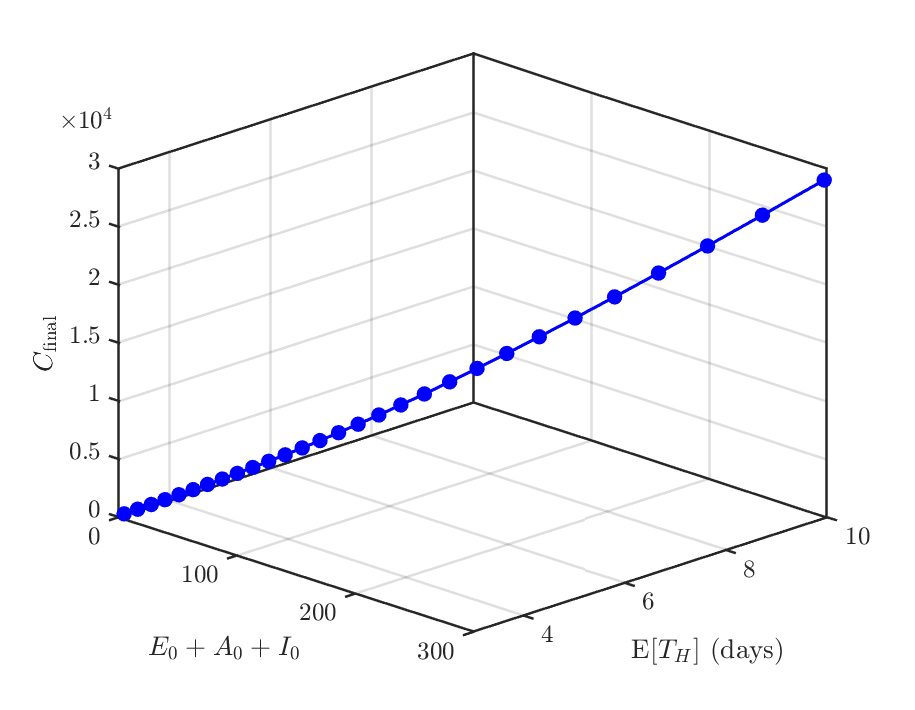}
	\caption{Effects of the IDH and the infection size at first detection on the final infection size.}
	\label{fig:size}
\end{figure}

The IDH affects the infection size at first outbreak detection and further influences the final infection size. Fig.~\ref{fig:size} illustrates the relationship among these three quantities. As the IDH increases, the duration of free transmission becomes longer, resulting in a larger number of infected individuals at the time of outbreak detection. These infected individuals then constitute the initial conditions of the post-detection intervention transmission model and influence the subsequent transmission process. Therefore, a longer IDH generally corresponds to a larger infection size at detection and a larger final infection size, indicating that earlier outbreak detection can reduce the accumulation of infected individuals during the free-transmission stage and limit subsequent transmission.

\section{Conclusions and discussion}\label{sec6}
This study proposes the concept of the IDH and develops a mathematical framework for calculating its distribution and mean. Based on an SEAIR free-transmission model, we derive key transmission characteristics and incidence using the renewal equation. By incorporating large-scale screening and medical diagnosis, we characterize the IDH and infection size at first detection. A post-detection intervention model is further developed to investigate subsequent transmission and the final infection size. Finally, the proposed modeling framework is applied to the 2021--2022 COVID-19 outbreak in Xi'an.

For the Xi'an outbreak, the mean IDH was approximately 6.6 days, with approximately 39 infected individuals remaining in states $E$, $A$, and $I$ at the mean detection time.  The results show that the IDH is jointly determined by the transmission process and the case-detection process. Increasing the frequency of nucleic acid screening or the probability that symptomatic infected individuals seek medical care can both shorten the mean IDH. When both medical diagnosis and large-scale screening are considered, outbreak detection relies mainly on medical diagnosis of symptomatic infected individuals once the screening interval reaches two weeks or longer. The infection age of the initial case at introduction has a non-monotonic effect on the mean IDH. In addition, within the parameter ranges considered, a larger basic reproduction number or a shorter mean generation interval corresponds to a shorter mean IDH. A longer IDH allows more infections to accumulate before detection, resulting in a larger final infection size and highlighting the importance of early outbreak detection in limiting the spread of emerging infectious diseases.

This study still has several limitations. First, the intervention model does not incorporate other control measures, such as contact tracing and preventive quarantine. Second, the transmission characteristics and case-detection processes may differ across emerging infectious diseases. Future work could incorporate additional interventions and disease-specific transmission and detection characteristics to further extend the proposed framework.

\bmhead{Acknowledgements}

This work was supported by the National Key Research and Development Program of China (No. 2025YFA1016500) and National Natural Science Foundation of China (Nos. 12231012, U23A20331).

\bmhead{Data availability}

The case data are obtained from the publicly available records of the Shaanxi Provincial Health Commission: \url{https://sxwjw.shaanxi.gov.cn/gab/wjyw/index_167.html}.

\section*{Declarations}
The authors declare no conflict of interest.

\begin{appendices}
	\section{Derivation of Eq.~\eqref{eq:GBnew}}
	\label{app:medical-diagnosis}
	
	We provide a detailed derivation of
	Eq.~\eqref{eq:GBnew}. For any $s\in(0,t]$,
	$\lambda(s)$ represents the incidence at time $s$.
	Assume that new infection events follow a
	nonhomogeneous Poisson process with intensity
	$\lambda(s)$ and that the medical diagnosis
	processes of different infected individuals
	are mutually independent.
	
	Divide the infection-age interval $(0,t]$
	into $n$ equal subintervals $0=t_0<t_1<\cdots<t_n=t,
	\
	\Delta t=\frac{t}{n}$. 
	In the $r$th subinterval $(t_{r-1},t_r]$, choose an arbitrary point $t_r^*\in(t_{r-1},t_r],\
	r=1,2,\ldots,n$. 
	The point $t_r^*$ represents an infection age
	within the $r$th subinterval. Individuals whose infection ages at time $t$
	lie in $(t_{r-1},t_r]$ were infected during
	the corresponding infection-time interval $
	[t-t_r,t-t_{r-1})$. 
	
	Let $N_r$ denote the number of new infections
	occurring in this infection-time interval.
	By the properties of the nonhomogeneous
	Poisson process,
	\[
	N_r\sim\operatorname{Poisson}(\mu_r),
	\]
	where
	\[
	\begin{aligned}
		\mu_r
		&=
		\int_{t-t_r}^{t-t_{r-1}}
		\lambda(s)\,ds=
		\int_{t_{r-1}}^{t_r}
		\lambda(t-a)\,da.
	\end{aligned}
	\]
	When $\Delta t$ is sufficiently small,
	we have $
	\mu_r
	\approx
	\lambda(t-t_r^*)\Delta t$.
	
	For individuals whose infection ages lie in the
	$r$th subinterval, the probability of remaining
	undetected through medical diagnosis by time $t$
	is approximated by $1-P_{T_B}(t_r^*)$. If $N_r=k$, the probability that all $k$
	individuals remain undetected through medical
	diagnosis by time $t$ is therefore approximated by
	$\left[
	1-P_{T_B}(t_r^*)
	\right]^k$.
	
	Averaging over the possible numbers of new
	infections, we obtain the following approximation
	to the probability that none of the individuals
	infected in the $r$th infection-time interval
	has been detected through medical diagnosis
	by time $t$
	\[
	\begin{aligned}
		G_{B,r}^{(n)}(t)
		&=
		\sum_{k=0}^{\infty}
		\Pr(N_r=k)
		\left[
		1-P_{T_B}(t_r^*)
		\right]^k
		\\
		&=
		\sum_{k=0}^{\infty}
		e^{-\mu_r}
		\frac{\mu_r^k}{k!}
		\left[
		1-P_{T_B}(t_r^*)
		\right]^k
		\\
		&=
		\exp\left\{
		-\mu_r P_{T_B}(t_r^*)
		\right\}.
	\end{aligned}
	\]
	Since the nonhomogeneous Poisson process has
	independent increments over disjoint time intervals
	and the medical diagnosis processes of different
	infected individuals are mutually independent,
	the corresponding approximation to the probability
	that none of the individuals infected during
	$(0,t]$ has been detected through medical
	diagnosis by time $t$ is
	\begin{equation}
		\begin{aligned}
			G_{B,\mathrm{new}}^{(n)}(t)
			&=
			\prod_{r=1}^{n}
			G_{B,r}^{(n)}(t)
			\\
			&=
			\prod_{r=1}^{n}
			\exp\left\{
			-\mu_r P_{T_B}(t_r^*)
			\right\}
			\\
			&=
			\exp\left\{
			-\sum_{r=1}^{n}
			\mu_r P_{T_B}(t_r^*)
			\right\}.
		\end{aligned}
		\label{eq:app-GB-partition}
	\end{equation}
	As $n\to\infty$, we have $\Delta t\to0$. Since $\mu_r
	=
	\int_{t_{r-1}}^{t_r}
	\lambda(t-a)\,da$, 
	taking the limit of the Riemann sum yields
	\[
	\begin{aligned}
		\lim_{n\to\infty}
		\sum_{r=1}^{n}
		\mu_r P_{T_B}(t_r^*)
		&=
		\lim_{n\to\infty}
		\sum_{r=1}^{n}
		P_{T_B}(t_r^*)
		\int_{t_{r-1}}^{t_r}
		\lambda(t-a)\,da=
		\int_0^t
		\lambda(t-a)P_{T_B}(a)\,da.
	\end{aligned}
	\]
	Consequently, taking the limit $n\to\infty$ in
	Eq.~\eqref{eq:app-GB-partition}, we obtain
	\[
	\begin{aligned}
		G_{B,\mathrm{new}}(t)
		&=
		\lim_{n\to\infty}
		G_{B,\mathrm{new}}^{(n)}(t)=
		\exp\left\{
		-\int_0^t
		\lambda(t-a)P_{T_B}(a)\,da
		\right\}.
	\end{aligned}
	\]
	This gives the probability that none of the individuals
	infected during $(0,t]$ has been detected through
	medical diagnosis by time $t$, as stated in
	Eq.~\eqref{eq:GBnew}.

\section{Intervention transmission model under exponential distributions }\label{secA2}
	
We now derive the intervention transmission model for the exponential-distribution case based on model~\eqref{eq:model_general}.

Substituting the exponential survival function of $\overline F_{T_E}(a)=e^{-\delta a}$and Eq.~\eqref{Qa} into Eq.~\eqref{eq:piE}, the probability that an individual infected after outbreak detection and the implementation of intervention measures remains in the latent state \(E\) at infection age \(a\) is
	\[
	\Pi_E(a)
	=
	e^{-(\delta+\kappa_E)a},
	\qquad
	\Pi_E(0)=1.
	\]
	Therefore,
	\[
	\frac{d\Pi_E(a)}{da}
	=
	-(\delta+\kappa_E)\Pi_E(a).
	\]
	Individuals in the latent state \(E\) at time \(t\) arise from two sources: individuals who are already in state \(E\) at outbreak detection and remain in state \(E\) at time \(t\), and individuals infected during \((0,t]\) who are still in state \(E\) at time \(t\). Therefore, the number of infected individuals in the latent state \(E\) at time \(t\) is
	\begin{equation}\label{Ae}
		\begin{aligned}
			E(t)
			&=
			E_0\Pi_E(t)+ \int_0^t
			\lambda(t-a)\Pi_E(a)\,da \\
			&=
			E_0 e^{-(\delta+\kappa_E)t}
			+
			\int_0^t
			\lambda(t-a)e^{-(\delta+\kappa_E)a}\,da.
		\end{aligned}
	\end{equation}
	Differentiating Eq.~\eqref{Ae} with respect to \(t\) gives
	\begin{equation}
		\begin{aligned}
			\frac{dE(t)}{dt}=
			\lambda(t)-(\delta+\kappa_E)E(t).
		\end{aligned}
		\label{Ade}
	\end{equation}
	
	Similarly, using the exponential survival function
	$\overline F_{T_E}(a)$ and Eq.~\eqref{Qa},
	the probability in Eq.~\eqref{eq:piEq} that an individual of infection age $a$ is in the isolated latent state $E_q$ becomes
	\begin{equation}\label{APiEQq}
		\Pi_{E_q}(a)
		=
		e^{-\delta a}\left(1-e^{-\kappa_Ea}\right).
	\end{equation}
	Moreover,
	\[
	\begin{aligned}
		\frac{d\Pi_{E_q}(a)}{da}
		&=
		\kappa_E\Pi_E(a)
		-
		\delta\Pi_{E_q}(a),
		\qquad
		\Pi_{E_q}(0)=0.
	\end{aligned}
	\]
	The number of infected individuals in the isolated latent state $E_q$ at time $t$ is
	\begin{equation}\label{Aeq}
		\begin{aligned}
			E_q(t)
			&=
			E_0\Pi_{E_q}(t)
			+
			\int_0^t
			\lambda(t-a)\Pi_{E_q}(a)\,da \\
			&=
			E_0e^{-\delta t}
			\left(1-e^{-\kappa_Et}\right)
			+
			\int_0^t
			\lambda(t-a)e^{-\delta a}
			\left(1-e^{-\kappa_Ea}\right)\,da.
		\end{aligned}
	\end{equation}
	Differentiating $E_q(t)$ with respect to $t$ yields
	\begin{equation}
		\begin{aligned}
			\frac{dE_q(t)}{dt}
			=
			\kappa_E E(t)-\delta E_q(t).
		\end{aligned}
		\label{Adeq}
	\end{equation}
	
	The probability that an individual of infection age $a$ is in the state $A$ is
	\begin{equation}
		\label{Aeq:PiA}
		\begin{aligned}
			\Pi_A(a)
			={}&
			\int_0^a
			\delta e^{-(\delta+\kappa_E)l}
			\Bigl[
			\rho e^{-(\eta_1+\kappa_A)(a-l)}
			+
			(1-\rho)e^{-(\eta_2+\kappa_A)(a-l)}
			\Bigr]
			\,dl .
		\end{aligned}
	\end{equation}
	Using the Leibniz rule for differentiation of integrals with a variable upper limit,
	\[
	\frac{d\Pi_A(a)}{da}
	=
	\delta\Pi_E(a)
	-\kappa_A\Pi_A(a)
	-
	\int_0^a
	\delta e^{-(\delta+\kappa_E)l}
	\Bigl[
	\rho\eta_1e^{-(\eta_1+\kappa_A)(a-l)}
	+
	(1-\rho)\eta_2e^{-(\eta_2+\kappa_A)(a-l)}
	\Bigr]
	\,dl,
	\]
	with
	\[
	\Pi_A(0)=0.
	\]
	For individuals who are already in state $A$ at the initial time,
	\[
	P_{AA}(t)
	=
	\rho e^{-(\eta_1+\kappa_A)t}
	+
	(1-\rho)e^{-(\eta_2+\kappa_A)t},
	\]
	and hence
	\[
	\frac{dP_{AA}(t)}{dt}
	=
	-\kappa_A P_{AA}(t)
	-
	\left[
	\rho\eta_1e^{-(\eta_1+\kappa_A)t}
	+
	(1-\rho)\eta_2e^{-(\eta_2+\kappa_A)t}
	\right].
	\]
	The number of infected individuals in the state $A$ at time $t$ is
	\[
	\begin{aligned}
		A(t)
		=
		E_0\Pi_A(t)
		+
		A_0P_{AA}(t)
		+
		\int_0^t
		\lambda(t-a)\Pi_A(a)\,da.
	\end{aligned}
	\]
	Differentiating $A(t)$ with respect to $t$ gives
	\begin{equation}
		\label{Aeq:dA_integral}
		\begin{aligned}
			\frac{dA(t)}{dt}
			&=
			\delta E(t)-\kappa_A A(t) \\
			&\quad
			-E_0
			\int_0^t
			\delta e^{-(\delta+\kappa_E)l}
			\Bigl[
			\rho\eta_1e^{-(\eta_1+\kappa_A)(t-l)}
			+
			(1-\rho)\eta_2e^{-(\eta_2+\kappa_A)(t-l)}
			\Bigr]\,dl \\
			&\quad
			-A_0
			\Bigl[
			\rho\eta_1e^{-(\eta_1+\kappa_A)t}
			+
			(1-\rho)\eta_2e^{-(\eta_2+\kappa_A)t}
			\Bigr] \\
			&\quad
			-\int_0^t
			\lambda(t-a)
			\int_0^a
			\delta e^{-(\delta+\kappa_E)l}
			\Bigl[
			\rho\eta_1e^{-(\eta_1+\kappa_A)(a-l)}
			+
			(1-\rho)\eta_2e^{-(\eta_2+\kappa_A)(a-l)}
			\Bigr]
			\,dl\,da .
		\end{aligned}
	\end{equation}
	
	The probability that an individual of infection age $a$ is in state $A_q$ is
	\begin{equation}
		\label{Aeq:PiAq}
		\begin{aligned}
			\Pi_{A_q}(a)
			={}&
			\int_0^a
			\delta e^{-\delta l}
			\Bigl[
			\rho e^{-\eta_1(a-l)}
			+
			(1-\rho)e^{-\eta_2(a-l)}
			\Bigr]\times
			\Bigl[
			1-e^{-\kappa_E l-\kappa_A(a-l)}
			\Bigr]
			\,dl .
		\end{aligned}
	\end{equation}
	This probability satisfies
	\[
	\begin{aligned}
		\frac{d\Pi_{A_q}(a)}{da}
		={}&
		\delta\Pi_{E_q}(a)
		+
		\kappa_A\Pi_A(a)
		\\
		&-
		\int_0^a
		\delta e^{-\delta l}
		\Bigl[
		\rho\eta_1e^{-\eta_1(a-l)}
		+
		(1-\rho)\eta_2e^{-\eta_2(a-l)}
		\Bigr]\times
		\Bigl[
		1-e^{-\kappa_E l-\kappa_A(a-l)}
		\Bigr]
		\,dl .
	\end{aligned}
	\]
	The number of infected individuals in the isolated state $A_q$ at time $t$ is
	\[
	\begin{aligned}
		A_q(t)
		&=
		A_q^{(0)}(t)+A_q^{\mathrm{new}}(t)=
		E_0\Pi_{A_q}(t)
		+
		A_0P_{AA_q}(t)
		+
		\int_0^t
		\lambda(t-a)\Pi_{A_q}(a)\,da.
	\end{aligned}
	\]
	Differentiating $A_q(t)$ with respect to $t$ gives
	\begin{equation}
		\label{Aeq:dAq_integral}
		\begin{aligned}
			\frac{dA_q(t)}{dt}
			&=
			\delta E_q(t)+\kappa_A A(t) \\
			&\quad
			-E_0
			\int_0^t
			\delta e^{-\delta l}
			\Bigl[
			\rho\eta_1e^{-\eta_1(t-l)}
			+
			(1-\rho)\eta_2e^{-\eta_2(t-l)}
			\Bigr]
			\\
			&\qquad\times
			\Bigl[
			1-e^{-\kappa_E l-\kappa_A(t-l)}
			\Bigr]\,dl \\
			&\quad
			-A_0
			\Bigl[
			\rho\eta_1e^{-\eta_1t}
			+
			(1-\rho)\eta_2e^{-\eta_2t}
			\Bigr]
			\Bigl[
			1-e^{-\kappa_At}
			\Bigr] \\
			&\quad
			-\int_0^t
			\lambda(t-a)
			\int_0^a
			\delta e^{-\delta l}
			\Bigl[
			\rho\eta_1e^{-\eta_1(a-l)}
			+
			(1-\rho)\eta_2e^{-\eta_2(a-l)}
			\Bigr]
			\\
			&\qquad\times
			\Bigl[
			1-e^{-\kappa_E l-\kappa_A(a-l)}
			\Bigr]
			\,dl\,da .
		\end{aligned}
	\end{equation}
	
	We next consider the probability that an individual of infection age $a$ is in
	the symptomatic infectious state $I$
	\begin{equation}
		\label{Aeq:PiI}
		\begin{aligned}
			\Pi_I(a)
			={}&
			\rho
			\int_0^a
			\int_0^{a-l}
			\delta e^{-(\delta+\kappa_E)l}
			\eta_1e^{-(\eta_1+\kappa_A)m}
			e^{-\gamma(a-l-m)}
			\\
			&\qquad\times
			\left[
			(1-\rho_b)
			+\rho_b e^{-\omega(a-l-m)}
			\right]
			\,dm\,dl.
		\end{aligned}
	\end{equation}
	Here, $\rho$ denotes the probability that an infected individual develops symptoms, and $\rho_b$ denotes the probability that a symptomatic infected individual seeks medical care and is diagnosed. In addition,
	\[
	\Pi_I(0)=0.
	\]
	Using the Leibniz rule,
	\begin{equation}
		\label{Aeq:dPiI}
		\begin{aligned}
			\frac{d\Pi_I(a)}{da}
			={}&
			\rho
			\int_0^a
			\delta e^{-(\delta+\kappa_E)l}
			\eta_1e^{-(\eta_1+\kappa_A)(a-l)}
			\,dl
			-\gamma\Pi_I(a)
			\\
			&-
			\rho\rho_b\omega
			\int_0^a
			\int_0^{a-l}
			\delta e^{-(\delta+\kappa_E)l}
			\eta_1e^{-(\eta_1+\kappa_A)m}
			e^{-(\gamma+\omega)(a-l-m)}
			\,dm\,dl .
		\end{aligned}
	\end{equation}
	For individuals who are in state $A$ at the initial time, the probability of being in the symptomatic infectious state $I$ at time $t$ is
	\[
	\begin{aligned}
		P_{AI}(t)
		={}&
		\rho
		\int_0^t
		\eta_1e^{-(\eta_1+\kappa_A)m}
		e^{-\gamma(t-m)}
		\left[
		(1-\rho_b)
		+
		\rho_b e^{-\omega(t-m)}
		\right]
		\,dm.
	\end{aligned}
	\]
	Differentiating $P_{AI}(t)$ with respect to $t$ gives
	\[
	\begin{aligned}
		\frac{dP_{AI}(t)}{dt}
		={}&
		\rho\eta_1e^{-(\eta_1+\kappa_A)t}
		-\gamma P_{AI}(t)-
		\rho\rho_b\omega
		\int_0^t
		\eta_1e^{-(\eta_1+\kappa_A)m}
		e^{-(\gamma+\omega)(t-m)}
		\,dm.
	\end{aligned}
	\]
	For individuals who are in state $I$ at the initial time, the probability of
	remaining in the symptomatic infectious state $I$ at time $t$ is
	\[
	P_{II}(t)
	=
	e^{-\gamma t}
	\left[
	(1-\rho_b)+\rho_b e^{-\omega t}
	\right],
	\]
	and
	\[
	\frac{dP_{II}(t)}{dt}
	=
	-\gamma P_{II}(t)
	-
	\rho_b\omega e^{-(\gamma+\omega)t}.
	\]
	Therefore, the number of symptomatic infectious individuals in state $I$ at time $t$
	is
	\[
	\begin{aligned}
		I(t)&=
		E_0\Pi_I(t)
		+
		A_0P_{AI}(t)
		+
		I_0P_{II}(t)
		+
		\int_0^t
		\lambda(t-a)\Pi_I(a)\,da.
	\end{aligned}
	\]
	Differentiating $I(t)$ with respect to $t$ gives
	\begin{equation}
		\label{Aeq:I_new}
		\begin{aligned}
			\frac{dI(t)}{dt}
			&=
			\rho\eta_1
			\Bigg[
			E_0
			\int_0^t
			\delta e^{-(\delta+\kappa_E)l}
			e^{-(\eta_1+\kappa_A)(t-l)}
			\,dl
			+
			A_0e^{-(\eta_1+\kappa_A)t}
			\\
			&\qquad+
			\int_0^t
			\lambda(t-a)
			\int_0^a
			\delta e^{-(\delta+\kappa_E)l}
			e^{-(\eta_1+\kappa_A)(a-l)}
			\,dl\,da
			\Bigg]
			\\
			&\quad
			-\rho\rho_b\omega\eta_1
			\Bigg[
			E_0
			\int_0^t
			\int_0^{t-l}
			\delta e^{-(\delta+\kappa_E)l}
			e^{-(\eta_1+\kappa_A)m}
			e^{-(\gamma+\omega)(t-l-m)}
			\,dm\,dl
			\\
			&\qquad+
			A_0
			\int_0^t
			e^{-(\eta_1+\kappa_A)m}
			e^{-(\gamma+\omega)(t-m)}
			\,dm
			\\
			&\qquad+
			\int_0^t
			\lambda(t-a)
			\int_0^a
			\int_0^{a-l}
			\delta e^{-(\delta+\kappa_E)l}
			e^{-(\eta_1+\kappa_A)m}
			e^{-(\gamma+\omega)(a-l-m)}
			\,dm\,dl\,da
			\Bigg]
			\\
			&\quad
			-\rho_b\omega I_0e^{-(\gamma+\omega)t}
			-\gamma I(t).
		\end{aligned}
	\end{equation}
	
	Individuals in state \(I_q\) arise from two sources. The first consists of infected individuals who are detected through nucleic acid testing before symptom onset and enter state \(A_q\), subsequently develop symptoms during the natural disease course, and then enter state \(I_q\). The second consists of individuals who are not detected through nucleic acid testing while in states \(E\) and \(A\), enter the symptomatic state \(I\), and are subsequently diagnosed through medical care and enter state \(I_q\). For an individual of infection age \(a\), the probability of being in the isolated symptomatic infectious state \(I_q\) is
	\begin{equation}
		\label{Aeq:PiIq}
		\begin{aligned}
			\Pi_{I_q}(a)
			={}&
			\rho
			\int_0^a
			\int_0^{a-l}
			\delta e^{-\delta l}
			\eta_1e^{-\eta_1m}
			e^{-\gamma(a-l-m)}
			\\
			&\qquad\times
			\Bigl\{
			1-e^{-\kappa_E l-\kappa_A m}
			+
			\rho_b e^{-\kappa_E l-\kappa_A m}
			\left[
			1-e^{-\omega(a-l-m)}
			\right]
			\Bigr\}
			\,dm\,dl.
		\end{aligned}
	\end{equation}
	Moreover,
	\[
	\Pi_{I_q}(0)=0.
	\]
	Differentiating $\Pi_{I_q}(a)$ with respect to infection age $a$ and applying
	the Leibniz rule gives
	\begin{equation}
		\label{Aeq:dPiIq}
		\begin{aligned}
			\frac{d\Pi_{I_q}(a)}{da}
			={}&
			\rho\eta_1
			\int_0^a
			\delta e^{-\delta l}
			e^{-\eta_1(a-l)}
			\left[
			1-e^{-\kappa_E l-\kappa_A(a-l)}
			\right]
			\,dl
			\\
			&+
			\rho\rho_b\omega
			\int_0^a
			\int_0^{a-l}
			\delta e^{-(\delta+\kappa_E)l}
			\eta_1e^{-(\eta_1+\kappa_A)m}
			e^{-(\gamma+\omega)(a-l-m)}
			\,dm\,dl
			\\
			&-
			\gamma\Pi_{I_q}(a).
		\end{aligned}
	\end{equation}
	For individuals who are in state $A$ at the initial time, the probability of being in the isolated symptomatic infectious state $I_q$ at time $t$ is
	\[
	\begin{aligned}
		P_{AI_q}(t)
		={}&
		\rho
		\int_0^t
		\eta_1e^{-\eta_1m}
		e^{-\gamma(t-m)}
		\left\{
		1-e^{-\kappa_A m}
		+
		\rho_b e^{-\kappa_A m}
		\left[
		1-e^{-\omega(t-m)}
		\right]
		\right\}
		\,dm.
	\end{aligned}
	\]
	Differentiating $P_{AI_q}(t)$ with respect to $t$ gives
	\[
	\begin{aligned}
		\frac{dP_{AI_q}(t)}{dt}
		={}&
		\rho\eta_1e^{-\eta_1t}
		\left(
		1-e^{-\kappa_A t}
		\right)
		-\gamma P_{AI_q}(t)
		\\
		&+
		\rho\rho_b\omega
		\int_0^t
		\eta_1e^{-(\eta_1+\kappa_A)m}
		e^{-(\gamma+\omega)(t-m)}
		\,dm.
	\end{aligned}
	\]
	For individuals who are in state $I$ at the initial time, the probability of being
	in state $I_q$ at time $t$ is
	\[
	P_{II_q}(t)
	=
	\rho_b e^{-\gamma t}
	\left[
	1-e^{-\omega t}
	\right],
	\]
	and
	\[
	\frac{dP_{II_q}(t)}{dt}
	=
	\rho_b\omega e^{-(\gamma+\omega)t}
	-\gamma P_{II_q}(t).
	\]
	The number of infected individuals in the isolated symptomatic infectious state $I_q$
	at time $t$ is
	\[
	\begin{aligned}
		I_q(t)&=
		E_0\Pi_{I_q}(t)
		+
		A_0P_{AI_q}(t)
		+
		I_0P_{II_q}(t)
		+
		\int_0^t
		\lambda(t-a)\Pi_{I_q}(a)\,da.
	\end{aligned}
	\]
	Differentiating $I_q(t)$ with respect to $t$ gives
	\begin{equation}
		\label{Aeq:I_q}
		\begin{aligned}
			\frac{dI_q(t)}{dt}&=
			\rho\eta_1
			\Bigg[
			E_0
			\int_0^t
			\delta e^{-\delta l}
			e^{-\eta_1(t-l)}
			\left(
			1-e^{-\kappa_E l-\kappa_A(t-l)}
			\right)
			\,dl
			\\
			&\qquad+
			A_0e^{-\eta_1t}
			\left(
			1-e^{-\kappa_A t}
			\right)
			\\
			&\qquad+
			\int_0^t
			\lambda(t-a)
			\int_0^a
			\delta e^{-\delta l}
			e^{-\eta_1(a-l)}
			\left(
			1-e^{-\kappa_E l-\kappa_A(a-l)}
			\right)
			\,dl\,da
			\Bigg]
			\\
			&+
			\rho\rho_b\omega\eta_1
			\Bigg[
			E_0
			\int_0^t
			\int_0^{t-l}
			\delta e^{-(\delta+\kappa_E)l}
			e^{-(\eta_1+\kappa_A)m}
			e^{-(\gamma+\omega)(t-l-m)}
			\,dm\,dl
			\\
			&\qquad+
			A_0
			\int_0^t
			e^{-(\eta_1+\kappa_A)m}
			e^{-(\gamma+\omega)(t-m)}
			\,dm
			\\
			&\qquad+
			\int_0^t
			\lambda(t-a)
			\int_0^a
			\int_0^{a-l}
			\delta e^{-(\delta+\kappa_E)l}
			e^{-(\eta_1+\kappa_A)m}
			e^{-(\gamma+\omega)(a-l-m)}
			\,dm\,dl\,da
			\Bigg]
			\\
			&+
			\rho_b\omega I_0
			e^{-(\gamma+\omega)t}
			-\gamma I_q(t).
		\end{aligned}
	\end{equation}
	
	The probability that an individual of infection age $a$ is in the recovered state $R$ is
	\begin{equation}
		\label{Aeq:PiR}
		\begin{aligned}
			\Pi_R(a)
			={}&
			\rho
			\int_0^a
			\int_0^{a-l}
			\int_0^{a-l-m}
			\delta e^{-\delta l}
			\eta_1e^{-\eta_1m}
			\gamma e^{-\gamma r}
			\,dr\,dm\,dl
			\\
			&+
			(1-\rho)
			\int_0^a
			\int_0^{a-l}
			\delta e^{-\delta l}
			\eta_2e^{-\eta_2m}
			\,dm\,dl .
		\end{aligned}
	\end{equation}
	Moreover,
	\[
	\Pi_R(0)=0.
	\]
	Differentiating $\Pi_R(a)$ with respect to $a$ gives
	\[
	\begin{aligned}
		\frac{d\Pi_R(a)}{da}
		={}&
		\rho
		\int_0^a
		\int_0^{a-l}
		\delta e^{-\delta l}
		\eta_1e^{-\eta_1m}
		\gamma e^{-\gamma(a-l-m)}
		\,dm\,dl
		\\
		&+
		(1-\rho)
		\int_0^a
		\delta e^{-\delta l}
		\eta_2e^{-\eta_2(a-l)}
		\,dl
		\\
		={}&
		\gamma
		\left[
		\Pi_I(a)+\Pi_{I_q}(a)
		\right]
		+
		(1-\rho)\eta_2
		\int_0^a
		\delta e^{-\delta l}
		e^{-\eta_2(a-l)}
		\,dl .
	\end{aligned}
	\]
	For individuals who are in state $A$ at the initial time, the probability of being in the recovered state $R$ at time $t$ is
	\[
	\begin{aligned}
		P_{AR}(t)
		={}&
		\rho
		\int_0^t
		\eta_1e^{-\eta_1m}
		\left[
		1-e^{-\gamma(t-m)}
		\right]
		\,dm+
		(1-\rho)
		\left[
		1-e^{-\eta_2t}
		\right].
	\end{aligned}
	\]
	Differentiating $P_{AR}(t)$ with respect to $t$ gives
	\[
	\begin{aligned}
		\frac{dP_{AR}(t)}{dt}
		={}&
		\rho\gamma
		\int_0^t
		\eta_1e^{-\eta_1m}
		e^{-\gamma(t-m)}
		\,dm
		+
		(1-\rho)\eta_2e^{-\eta_2t}.
	\end{aligned}
	\]
	For individuals who are in state $I$ at the initial time, the probability of being in the recovered state $R$ at time $t$ is
	\[
	P_{IR}(t)
	=
	1-e^{-\gamma t},
	\]
	and
	\[
	\frac{dP_{IR}(t)}{dt}
	=
	\gamma e^{-\gamma t}.
	\]
	Therefore, the number of individuals in the recovered state $R$ at time $t$ is
	\[
	\begin{aligned}
		R(t)
		&=
		R(0)
		+
		E_0\Pi_R(t)
		+
		A_0P_{AR}(t)
		+
		I_0P_{IR}(t)
		+
		\int_0^t
		\lambda(t-a)\Pi_R(a)\,da.
	\end{aligned}
	\]
	Differentiating $R(t)$ with respect to $t$ gives
	\begin{equation}
		\label{Aeq:R_eq}
		\begin{aligned}
			\frac{dR(t)}{dt}
			&=
			\gamma
			\left[
			I(t)+I_q(t)
			\right]+
			(1-\rho)\eta_2
			\Bigg[
			E_0
			\int_0^t
			\delta e^{-\delta l}
			e^{-\eta_2(t-l)}
			\,dl
			+
			A_0e^{-\eta_2t}
			\\
			&\qquad+
			\int_0^t
			\lambda(t-a)
			\int_0^a
			\delta e^{-\delta l}
			e^{-\eta_2(a-l)}
			\,dl\,da
			\Bigg].
		\end{aligned}
	\end{equation}
	Therefore, when the durations of all disease stages follow exponential distributions,
	using Eqs.~(\ref{Ade}), (\ref{Adeq}), (\ref{Aeq:dA_integral}),
	(\ref{Aeq:dAq_integral}), (\ref{Aeq:I_new}), (\ref{Aeq:I_q}),
	and~(\ref{Aeq:R_eq}), the intervention transmission model incorporating isolation
	measures can be written as
	\begin{equation}
		\left\{
		\begin{aligned}
			\frac{dS(t)}{dt}
			={}&
			-\beta S(t)\left[A(t)+\xi_I I(t)\right],
			\\[1mm]
			\frac{dE(t)}{dt}
			={}&
			\beta S(t)\left[A(t)+\xi_I I(t)\right]
			-(\delta+\kappa_E)E(t),
			\\[1mm]
			\frac{dE_q(t)}{dt}
			={}&
			\kappa_EE(t)-\delta E_q(t),
			\\[1mm]
			\frac{dA(t)}{dt}
			={}&
			\delta E(t)-\kappa_AA(t)
			-\rho\eta_1 K_{1}(t)
			-(1-\rho)\eta_2 K_{2}(t),
			\\[1mm]
			\frac{dA_q(t)}{dt}
			={}&
			\delta E_q(t)+\kappa_AA(t)
			-\rho\eta_1 K_{q1}(t)
			-(1-\rho)\eta_2 K_{q2}(t),
			\\[1mm]
			\frac{dI(t)}{dt}
			={}&
			\rho\eta_1 K_{1}(t)
			-\rho_b\omega\mathcal{J}_{\omega}(t)
			-\gamma I(t),
			\\[1mm]
			\frac{dI_q(t)}{dt}
			={}&
			\rho\eta_1 K_{q1}(t)
			+\rho_b\omega\mathcal{J}_{\omega}(t)
			-\gamma I_q(t),
			\\[1mm]
			\frac{dR(t)}{dt}
			={}&
			\gamma\left[I(t)+I_q(t)\right]
			+
			(1-\rho)\eta_2
			\left[
			K_{2}(t)+K_{q2}(t)
			\right].
		\end{aligned}
		\right.
		\label{Aeq:jw}
	\end{equation}
	where
	\begin{equation}
		\begin{aligned}
			K_{1}(t)
			={}&
			E_0
			\int_{0}^{t}
			\delta e^{-(\delta+\kappa_E)l}
			e^{-(\eta_1+\kappa_A)(t-l)}
			\,dl
			+
			A_0e^{-(\eta_1+\kappa_A)t}
			\\
			&+
			\int_{0}^{t}
			\lambda(t-a)
			\int_{0}^{a}
			\delta e^{-(\delta+\kappa_E)l}
			e^{-(\eta_1+\kappa_A)(a-l)}
			\,dl\,da .
		\end{aligned}
		\label{Aeq:K1}
	\end{equation}
	\begin{equation}
		\begin{aligned}
			K_{2}(t)
			={}&
			E_0
			\int_{0}^{t}
			\delta e^{-(\delta+\kappa_E)l}
			e^{-(\eta_2+\kappa_A)(t-l)}
			\,dl
			+
			A_0e^{-(\eta_2+\kappa_A)t}
			\\
			&+
			\int_{0}^{t}
			\lambda(t-a)
			\int_{0}^{a}
			\delta e^{-(\delta+\kappa_E)l}
			e^{-(\eta_2+\kappa_A)(a-l)}
			\,dl\,da .
		\end{aligned}
		\label{Aeq:K2}
	\end{equation}
	\begin{equation}
		\begin{aligned}
			K_{q1}(t)
			={}&
			E_0
			\int_{0}^{t}
			\delta e^{-\delta l}
			e^{-\eta_1(t-l)}
			\left[
			1-e^{-\kappa_E l-\kappa_A(t-l)}
			\right]
			\,dl
			\\
			&+
			A_0e^{-\eta_1t}
			\left[
			1-e^{-\kappa_At}
			\right]
			\\
			&+
			\int_{0}^{t}
			\lambda(t-a)
			\int_{0}^{a}
			\delta e^{-\delta l}
			e^{-\eta_1(a-l)}
			\left[
			1-e^{-\kappa_E l-\kappa_A(a-l)}
			\right]
			\,dl\,da .
		\end{aligned}
		\label{Aeq:Kq1}
	\end{equation}
	\begin{equation}
		\begin{aligned}
			K_{q2}(t)
			={}&
			E_0
			\int_{0}^{t}
			\delta e^{-\delta l}
			e^{-\eta_2(t-l)}
			\left[
			1-e^{-\kappa_E l-\kappa_A(t-l)}
			\right]
			\,dl
			\\
			&+
			A_0e^{-\eta_2t}
			\left[
			1-e^{-\kappa_At}
			\right]
			\\
			&+
			\int_{0}^{t}
			\lambda(t-a)
			\int_{0}^{a}
			\delta e^{-\delta l}
			e^{-\eta_2(a-l)}
			\left[
			1-e^{-\kappa_E l-\kappa_A(a-l)}
			\right]
			\,dl\,da .
		\end{aligned}
		\label{Aeq:Kq2}
	\end{equation}
	\begin{equation}
		\begin{aligned}
			\mathcal{J}_{\omega}(t)
			={}&
			\rho E_0
			\int_{0}^{t}
			\int_{0}^{t-l}
			\delta e^{-(\delta+\kappa_E)l}
			\eta_1e^{-(\eta_1+\kappa_A)m}
			e^{-(\gamma+\omega)(t-l-m)}
			\,dm\,dl
			\\
			&+
			\rho A_0
			\int_{0}^{t}
			\eta_1e^{-(\eta_1+\kappa_A)m}
			e^{-(\gamma+\omega)(t-m)}
			\,dm
			+
			I_0e^{-(\gamma+\omega)t}
			\\
			&+
			\rho
			\int_{0}^{t}
			\lambda(t-a)
			\int_{0}^{a}
			\int_{0}^{a-l}
			\delta e^{-(\delta+\kappa_E)l}
			\eta_1e^{-(\eta_1+\kappa_A)m}
			e^{-(\gamma+\omega)(a-l-m)}
			\,dm\,dl\,da .
		\end{aligned}
		\label{Aeq:Jomega}
	\end{equation}
	Here, \(\rho K_1(t)\) represents the number of infected individuals who will eventually develop symptoms and are still in state \(A\) at time \(t\), whereas \((1-\rho)K_2(t)\) represents the number of infected individuals who remain asymptomatic throughout the infectious period and are still in state \(A\) at time \(t\). Similarly, \(\rho K_{q1}(t)\) and \((1-\rho)K_{q2}(t)\) represent the corresponding two groups of infected individuals in the isolated state \(A_q\). \(\rho_b\omega\mathcal{J}_{\omega}(t)\) represents the rate at which symptomatic infected individuals transition from state \(I\) to state \(I_q\) through medical diagnosis at time \(t\).
	
	We next derive the ordinary differential equations satisfied by the integral auxiliary functions $K_{1}(t)$, $K_{2}(t)$, $K_{q1}(t)$, $K_{q2}(t)$, and $\mathcal{J}_{\omega}(t)$. First, define the probability kernel corresponding to Eq.~\eqref{Aeq:K1} as
	\[
	\Pi_{A,1}(a)
	=
	\int_{0}^{a}
	\delta e^{-(\delta+\kappa_E)l}
	e^{-(\eta_1+\kappa_A)(a-l)}
	\,dl .
	\]
	Then Eq.~\eqref{Aeq:K1} can be written as
	\[
	\begin{aligned}
		K_{1}(t)
		={}&
		E_0\Pi_{A,1}(t)
		+
		A_0e^{-(\eta_1+\kappa_A)t}
		+
		\int_{0}^{t}
		\lambda(t-a)\Pi_{A,1}(a)
		\,da .
	\end{aligned}
	\]
	By the Leibniz rule,
	\[
	\frac{d\Pi_{A,1}(a)}{da}
	=
	\delta e^{-(\delta+\kappa_E)a}
	-
	(\eta_1+\kappa_A)\Pi_{A,1}(a),
	\qquad
	\Pi_{A,1}(0)=0.
	\]
	Hence,
	\[
	\begin{aligned}
		\frac{dK_{1}(t)}{dt}
		={}&
		\delta
		\left[
		E_0e^{-(\delta+\kappa_E)t}
		+
		\int_{0}^{t}
		\lambda(t-a)e^{-(\delta+\kappa_E)a}
		\,da
		\right]
		\\
		&-
		(\eta_1+\kappa_A)
		\left[
		E_0\Pi_{A,1}(t)
		+
		A_0e^{-(\eta_1+\kappa_A)t}
		+
		\int_{0}^{t}
		\lambda(t-a)\Pi_{A,1}(a)
		\,da
		\right].
	\end{aligned}
	\]
	Since
	\[
	E(t)
	=
	E_0e^{-(\delta+\kappa_E)t}
	+
	\int_{0}^{t}
	\lambda(t-a)e^{-(\delta+\kappa_E)a}
	\,da,
	\]
	we obtain
	\begin{equation}
		\frac{dK_{1}(t)}{dt}
		=
		\delta E(t)
		-
		(\eta_1+\kappa_A)K_{1}(t).
		\label{Aeq:dk1}
	\end{equation}
	
	Similarly, define the probability kernel corresponding to Eq.~\eqref{eq:K2} as
	\[
	\Pi_{A,2}(a)
	=
	\int_{0}^{a}
	\delta e^{-(\delta+\kappa_E)l}
	e^{-(\eta_2+\kappa_A)(a-l)}
	\,dl ,
	\]
	Then Eq.~\eqref{Aeq:K2} can be written as
	\[
	\begin{aligned}
		K_{2}(t)
		={}&
		E_0\Pi_{A,2}(t)
		+
		A_0e^{-(\eta_2+\kappa_A)t}+
		\int_{0}^{t}
		\lambda(t-a)\Pi_{A,2}(a)
		\,da .
	\end{aligned}
	\]
	Differentiating  $\Pi_{A,2}(a)$  with respect to $a$ gives
	\[
	\frac{d\Pi_{A,2}(a)}{da}
	=
	\delta e^{-(\delta+\kappa_E)a}
	-
	(\eta_2+\kappa_A)\Pi_{A,2}(a),
	\qquad
	\Pi_{A,2}(0)=0.
	\]
	Hence,
	\begin{equation}
		\frac{dK_{2}(t)}{dt}
		=
		\delta E(t)
		-
		(\eta_2+\kappa_A)K_{2}(t).
		\label{Aeq:dK2}
	\end{equation}

	Next, define the probability kernel corresponding to Eq.~\eqref{eq:Kq1} as
	\[
	\Pi_{A_q,1}(a)
	=
	\int_{0}^{a}
	\delta e^{-\delta l}
	e^{-\eta_1(a-l)}
	\left[
	1-e^{-\kappa_E l-\kappa_A(a-l)}
	\right]
	\,dl .
	\]
	Then Eq.~\eqref{Aeq:Kq1} can be written as
	\[
	\begin{aligned}
		K_{q1}(t)
		={}&
		E_0\Pi_{A_q,1}(t)
		+
		A_0e^{-\eta_1t}
		\left(
		1-e^{-\kappa_At}
		\right)
		+
		\int_{0}^{t}
		\lambda(t-a)\Pi_{A_q,1}(a)
		\,da .
	\end{aligned}
	\]
	Differentiating $\Pi_{A_q,1}(a)$ with respect to $a$ gives
	\[
	\begin{aligned}
		\frac{d\Pi_{A_q,1}(a)}{da}
		={}&
		\delta e^{-\delta a}
		\left(
		1-e^{-\kappa_Ea}
		\right)
		+
		\kappa_A\Pi_{A,1}(a)
		-
		\eta_1\Pi_{A_q,1}(a).
	\end{aligned}
	\]
	Thus,
	\[
	\begin{aligned}
		\frac{dK_{q1}(t)}{dt}
		={}&
		\delta
		\left[
		E_0e^{-\delta t}
		\left(
		1-e^{-\kappa_Et}
		\right)
		+
		\int_{0}^{t}
		\lambda(t-a)e^{-\delta a}
		\left(
		1-e^{-\kappa_Ea}
		\right)
		\,da
		\right]
		\\
		&+
		\kappa_A
		\left[
		E_0\Pi_{A,1}(t)
		+
		A_0e^{-(\eta_1+\kappa_A)t}
		+
		\int_{0}^{t}
		\lambda(t-a)\Pi_{A,1}(a)
		\,da
		\right]
		\\
		&-
		\eta_1
		\left[
		E_0\Pi_{A_q,1}(t)
		+
		A_0e^{-\eta_1t}
		\left(
		1-e^{-\kappa_At}
		\right)
		+
		\int_{0}^{t}
		\lambda(t-a)\Pi_{A_q,1}(a)
		\,da
		\right].
	\end{aligned}
	\]
	Since
	\[
	E_q(t)
	=
	E_0e^{-\delta t}
	\left(
	1-e^{-\kappa_Et}
	\right)
	+
	\int_{0}^{t}
	\lambda(t-a)e^{-\delta a}
	\left(
	1-e^{-\kappa_Ea}
	\right)
	\,da,
	\]
	we obtain
	\begin{equation}
		\frac{dK_{q1}(t)}{dt}
		=
		\delta E_q(t)
		+
		\kappa_AK_{1}(t)
		-
		\eta_1K_{q1}(t).
		\label{Aeq:dKq1}
	\end{equation}
	
	Similarly, denote the probability kernel corresponding to Eq.~\eqref{Aeq:Kq2} by
	\[
	\Pi_{A_q,2}(a)
	=
	\int_{0}^{a}
	\delta e^{-\delta l}
	e^{-\eta_2(a-l)}
	\left[
	1-e^{-\kappa_E l-\kappa_A(a-l)}
	\right]
	\,dl .
	\]
	Then,
	\[
	\begin{aligned}
		K_{q2}(t)
		={}&
		E_0\Pi_{A_q,2}(t)
		+
		A_0e^{-\eta_2t}
		\left(
		1-e^{-\kappa_At}
		\right)+
		\int_{0}^{t}
		\lambda(t-a)\Pi_{A_q,2}(a)
		\,da .
	\end{aligned}
	\]
	Differentiating \(\Pi_{A_q,2}(a)\) with respect to \(a\) gives
	\[
	\begin{aligned}
		\frac{d\Pi_{A_q,2}(a)}{da}
		={}&
		\delta e^{-\delta a}
		\left(
		1-e^{-\kappa_Ea}
		\right)
		+
		\kappa_A\Pi_{A,2}(a)
		-
		\eta_2\Pi_{A_q,2}(a),
		\\
		&\qquad
		\Pi_{A_q,2}(0)=0.
	\end{aligned}
	\]
	Therefore,
	\begin{equation}
		\frac{dK_{q2}(t)}{dt}
		=
		\delta E_q(t)
		+
		\kappa_AK_{2}(t)
		-
		\eta_2K_{q2}(t).
		\label{Aeq:dKq2}
	\end{equation}
	
	Finally, consider the auxiliary function \(\mathcal{J}_{\omega}(t)\). Define the probability kernel in Eq.~\eqref{Aeq:Jomega} corresponding to an individual initially entering state \(E\) as
	\[
	\Pi_{I,\omega}(a)
	=
	\int_{0}^{a}
	\int_{0}^{a-l}
	\delta e^{-(\delta+\kappa_E)l}
	\eta_1e^{-(\eta_1+\kappa_A)m}
	e^{-(\gamma+\omega)(a-l-m)}
	\,dm\,dl ,
	\]
	and define the probability kernel corresponding to an individual who is in state \(A\) at the initial time as
	\[
	P_{AI,\omega}(t)
	=
	\int_{0}^{t}
	\eta_1e^{-(\eta_1+\kappa_A)m}
	e^{-(\gamma+\omega)(t-m)}
	\,dm .
	\]
	Then Eq.~\eqref{Aeq:Jomega} can be written as
	\[
	\begin{aligned}
		\mathcal{J}_{\omega}(t)
		={}&
		\rho E_0\Pi_{I,\omega}(t)
		+
		\rho A_0P_{AI,\omega}(t)
		+
		I_0e^{-(\gamma+\omega)t}+
		\rho
		\int_{0}^{t}
		\lambda(t-a)\Pi_{I,\omega}(a)
		\,da .
	\end{aligned}
	\]
	Differentiating \(\Pi_{I,\omega}(a)\) with respect to \(a\), and applying Leibniz's rule for integrals with variable upper limits, gives
	\[
	\frac{d\Pi_{I,\omega}(a)}{da}
	=
	\eta_1\Pi_{A,1}(a)
	-
	(\gamma+\omega)\Pi_{I,\omega}(a),
	\qquad
	\Pi_{I,\omega}(0)=0.
	\]
	Similarly,
	\[
	\frac{dP_{AI,\omega}(t)}{dt}
	=
	\eta_1e^{-(\eta_1+\kappa_A)t}
	-
	(\gamma+\omega)P_{AI,\omega}(t),
	\qquad
	P_{AI,\omega}(0)=0.
	\]
	It therefore follows that
	\[
	\begin{aligned}
		\frac{d\mathcal{J}_{\omega}(t)}{dt}
		={}&
		\rho\eta_1
		\Bigg[
		E_0\Pi_{A,1}(t)
		+
		A_0e^{-(\eta_1+\kappa_A)t}
		+
		\int_{0}^{t}
		\lambda(t-a)\Pi_{A,1}(a)
		\,da
		\Bigg]
		\\
		&-
		(\gamma+\omega)
		\Bigg[
		\rho E_0\Pi_{I,\omega}(t)
		+
		\rho A_0P_{AI,\omega}(t)
		+
		I_0e^{-(\gamma+\omega)t}+
		\rho
		\int_{0}^{t}
		\lambda(t-a)\Pi_{I,\omega}(a)
		\,da
		\Bigg].
	\end{aligned}
	\]
	Using the definitions of \(K_1(t)\) and \(\mathcal{J}_{\omega}(t)\), we obtain
	\begin{equation}
		\frac{d\mathcal{J}_{\omega}(t)}{dt}
		=
		\rho\eta_1K_{1}(t)
		-
		(\gamma+\omega)\mathcal{J}_{\omega}(t).
		\label{Aeq:dJomega}
	\end{equation}

	In summary, under the above population-level averaging and exponential-equivalent approximation, combining equations~\eqref{Aeq:dk1}--\eqref{Aeq:dJomega} for the five
	integral auxiliary functions with model~\eqref{Aeq:jw}, the exponential-equivalent
	model~\eqref{Aeq:jw} can be further written as the following closed system of ordinary
	differential equations:
	\begin{equation}
		\left\{
		\begin{aligned}
			\frac{dS(t)}{dt}
			&=
			-\beta S(t)
			\left[
			\rho K_1(t)
			+
			(1-\rho)K_2(t)
			+
			\xi_I I(t)
			\right],
			\\ 
			\frac{dE(t)}{dt}
			&=
			\beta S(t)
			\left[
			\rho K_1(t)
			+
			(1-\rho)K_2(t)
			+
			\xi_I I(t)
			\right]
			-
			(\delta+\kappa_E)E(t),
			\\ 
			\frac{dE_q(t)}{dt}
			&=
			\kappa_EE(t)-\delta E_q(t),
			\\ 
			\frac{dK_1(t)}{dt}
			&=
			\delta E(t)
			-
			(\eta_1+\kappa_A)K_1(t),
			\\ 
			\frac{dK_2(t)}{dt}
			&=
			\delta E(t)
			-
			(\eta_2+\kappa_A)K_2(t),
			\\ 
			\frac{dK_{q1}(t)}{dt}
			&=
			\delta E_q(t)
			+
			\kappa_AK_1(t)
			-
			\eta_1K_{q1}(t),
			\\ 
			\frac{dK_{q2}(t)}{dt}
			&=
			\delta E_q(t)
			+
			\kappa_AK_2(t)
			-
			\eta_2K_{q2}(t),
			\\ 
			\frac{d\mathcal J_{\omega}(t)}{dt}
			&=
			\rho\eta_1K_1(t)
			-
			(\gamma+\omega)\mathcal J_{\omega}(t),
			\\ 
			\frac{dI(t)}{dt}
			&=
			\rho\eta_1K_1(t)
			-
			\rho_b\omega\mathcal J_{\omega}(t)
			-
			\gamma I(t),
			\\ 
			\frac{dI_q(t)}{dt}
			&=
			\rho\eta_1K_{q1}(t)
			+
			\rho_b\omega\mathcal J_{\omega}(t)
			-
			\gamma I_q(t),
			\\ 
			\frac{dR(t)}{dt}
			&=
			\gamma
			\left[
			I(t)+I_q(t)
			\right]
			+
			(1-\rho)\eta_2
			\left[
			K_2(t)+K_{q2}(t)
			\right].
		\end{aligned}
		\right.
		\label{Aeqode}
	\end{equation}
	The corresponding initial conditions are
	\[
	S(0)=S_0,\quad
	E(0)=E_0,\quad
	E_q(0)=0,
	\]
	\[
	K_1(0)=K_2(0)=A_0,\quad
	K_{q1}(0)=K_{q2}(0)=0,
	\]
	and
	\[
	\mathcal J_{\omega}(0)=I_0,\quad
	I(0)=I_0,\quad
	I_q(0)=0.
	\]
	Here,
	\[
	\kappa_E
	=
	-q\ln(1-\bar g_E),
	\quad
	\kappa_A
	=
	-q\ln(1-\bar g_A),
	\quad
	q=\frac{1}{\Delta T}.
	\]

\end{appendices}

\bibliography{references}%

\end{document}